\documentclass[final,3p,sort&compress]{elsarticle}

\usepackage{hyperref}
\usepackage{graphicx}
\usepackage{xcolor}
\usepackage{amsmath}
\usepackage{mathtools}
\usepackage{amsthm}
\usepackage{yhmath}
\usepackage{delimset}
\usepackage{empheq}
\usepackage{units}
\usepackage{amssymb}

\allowdisplaybreaks

\makeatletter
\def\NAT@spacechar{}
\makeatother

\usepackage{titlesec}
\titleformat*{\subsection}{\bfseries}

\journal{Xxxxx}

\begin{document}

\begin{frontmatter}

\title{Tempered Fractional Brownian Motion Revisited Via Fractional Ornstein-Uhlenbeck Processes 
}

\author[label1]{S.C. Lim\corref{cor1}}
\address[label1]{50 Holland Road,
\#02-01 Botanika,
Singapore 258853
}
\ead{sclim47@gmail.com}
\cortext[cor1]{corresponding author}

\author[label5]{Chai Hok Eab}
\address[label5]{
7/9 Rongmuang 5, Pathumwan, Bangkok 10330,Thailand
}
\ead{Chaihok.E@gmail.com}

\begin{abstract}
Tempered fractional Brownian motion is revisited from the viewpoint of reduced fractional Ornstein-Uhlenbeck process. 
Many of the basic properties of the tempered fractional Brownian motion can be shown to be direct consequences or modifications of the properties
of fractional Ornstein-Uhlenbeck process. 
Mixed tempered fractional Brownian motion is introduced and its properties are considered. 
Tempered fractional Brownian motion is generalised from single index to two indices. 
Finally, tempered multifractional Brownian motion and its properties are studied.
\end{abstract}


\end{frontmatter}


 
\section{Introduction}
\label{sec:TFBM_Intro}\noindent
Recently an interesting new stochastic process called tempered fractional Brownian motion (TFBM) has been introduced by Meerschaert and Sabzikar 
\cite{MeerschaertSbzikar2013,MeerschaertSabzikar2014}.
TFBM can be defined by modifying the moving average representation of a fractional Brownian motion (FBM) 
with the inclusion of an exponential tempering factor to the power-law kernel. 
This process has some nice features such as scaling property, stationary increments, moving average representation,  etc. 
It has an advantage in modelling realistic data as the tempering parameter that modifying the power law kernel can be chosen to 
the desired degree of accuracy over a finite interval.

In view of the increasing interest in TFBM, we would like to provide an alternative approach to 
TFBM.
It is possible to define TFBM as the reduced fractional Ornstein-Uhlenbeck process (RFOU), 
in the same sense as the modified ``reduced'' process introduced by Mandelbrot and van Ness \cite{MandelbrotVanNess68}
to get rid of divergence in the definition of fractional Brownian motion (FBM) based on Liouville-Weyl fractional integration of white noise. 
In contrast, the fractional Ornstein-Uhlenbeck process (FOU) is well-defined based on the Liouville-Weyl fractional integral. 
Instead, one encounters divergence in the ``massless'' limit of the covariance and variance of FOU 
in the attempt to recover FBM (see section 3 for more details). 
The main advantage of treating TFBM as RFOU is that many of its basic properties are inherited directly or modified from the properties of FOU,  
which have been well-studied 
\cite{Eab06a,Eab06b,LimLITeo2007,LimLiTeo2008,LimTeo2009c}.
In addition, such a setting allows the possible generalisation of TFBM from single index to two indices. 

In this paper three types of generalisations of TFBM are studied. 
Extension from TFBM to mixed TFBM can be treated in the  same way as mixed FBM 
\cite{Cheridito2001,El-Nouty2003}.
RFOU of single index can be generalised to two indices and the resulting process can be regarded as TFBM with two indices. 
Generalisation of TFBM to tempered multifractional Brownian motion (TMBM) with variable index can be carried out in the same manner as the extension of FBM to multifractional Brownian motion 
\cite{Levy-VehelPeltier1996,BenassiJaffardRoux1997}.

The outline of this paper is as follows. FOU and its properties are summarised in Section {\bf \ref{sec:fracOU}}. 
In the subsequent section, TFBM is defined in terms of RFOU, and most of its properties can be verified based on that of FOU. 
Mixed TFBM is discussed in Section {\bf \ref{sec:TFBMRFOU}}. 
Extension of TFBM to two indices is given in Section {\bf \ref{sec:MixTFBM}}. 
TFBM with variable index or TMFM is studied in Section {\bf\ref{sec:2indRFOU}}. 
The last section provides some concluding remarks.


\section{Fractional Ornstein-Uhlenbeck Process}
\label{sec:fracOU}
\noindent
There are several ways to generalise the Ornstein-Uhlenbeck process to its fractional counterpart 
\cite{CheriditoKawaguchiMaejima2003,LimMunuandy2003,Magdziarz2008}.
In this paper, 
the fractional generalisation of Ornstein-Uhlenbeck process will be considered as the solution of the following Langevin equation:
\begin{align}
  \bigl({_aD_t}+\lambda\bigr)^\alpha X_{\alpha,\lambda}(t) & = \eta(t), & 
         \alpha & > 1/2, 
\label{eq:fracOU_0010}
\end{align}
where $\lambda > 0$ is a positive constant, and $\eta(t)$ is the Gaussian white noise with zero mean and delta-correlated covariance. 
The condition $\alpha > 1/2$ is imposed to ensure the solution has finite variance. 
The fractional derivative ${_aD_t^\alpha}$ is defined by 
\cite{Samko1993,Kilbas2006}
\begin{align}
  {_aD_t^\alpha}f(t) & = \frac{1}{\Gamma(n-\alpha)} \left(\frac{d}{dt}\right)^n  \int_a^t \frac{f(u)}{(t-u)^{\alpha-n+1}} du, &
          & n-1 < \alpha < n. 
\label{eq:fracOU_0020}
\end{align}
For $a =0$, the fractional derivative is known as the Riemann-Liouville fractional derivative; 
when $a = -\infty$, 
it is called the Liouville-Weyl fractional derivative. 
Note that in some physics literature the Liouville-Weyl fractional derivative is also known as Weyl derivative, 
which is preferred in this paper.    

One can formally defined the ``shifted'' fractional derivative 
$ \bigl({_aD_t} + \lambda\bigr)^\alpha$ in terms of the unshifted derivative ${_aD_t^\alpha}$.
By using binomial expansion, it is possible to express the shifted fractional derivative in terms of unshifted ones:  
\begin{align}
  \bigl({_aD_t} + \lambda\bigr)^\alpha & = \sum_{j=0}^\infty \binom{\alpha}{j} \lambda^j {_aD_t^{\alpha-j}}.
\label{eq:fracOU_0030}
\end{align}
Note that this shifted fractional differential operator is closely related to the one-dimensional Bessel fractional
derivative and Bessel potential (Samko et al \cite{Samko1993}, page 336).
It can also be treated in a more
rigorous way by using hypersingular integrals \cite{Samko2002}.

The following operator identity holds for both the Riemann-Liouville and Weyl
fractional derivatives \cite{Eab06b}:
\begin{align}
   \bigl({_aD_t} + \lambda\bigr)^\alpha & = e^{-\lambda{t}} {_aD_t^{\alpha}} e^{\lambda{t}} .
\label{eq:fracOU_0040}
\end{align}
It can be verified by using the Binomial expansion and applying the generalized Leibniz rule, 
or by considering its Fourier transform in the case of Weyl fractional derivative, 
and Laplace transform for the Riemann-Liouville fractional derivative. 
Note that the shifted fractional derivative 
(\ref{eq:fracOU_0040})
is also known as tempered fractional derivative in subsequent work 
\cite[see for examples,\ ][]{CarteaCastillo-Negrete2007,SabzikarMeerschaertChen2015,CaoLiChen2014}.

Using 
(\ref{eq:fracOU_0040})
the fractional Langevin equation (\ref{eq:fracOU_0010}) can be re-expressed as
\begin{align}
    \bigl(e^{-\lambda{t}}{_aD_t^\alpha}e^{\lambda{t}}\bigr)X_{\alpha,\lambda}(t) & = \eta(t).
\label{eq:fracOU_0050}
\end{align}
The solution is given by
\begin{align}
  X_{\alpha,\lambda}(t) & = = \Bigl( e^{-\lambda{t}} {_aD_t^{\alpha}} e^{\lambda{t}} \Bigr)^{-1}\eta(t)
                       = \Bigl( e^{-\lambda{t}} {_aI_t^{\alpha}} e^{\lambda{t}} \Bigr)\eta(t) \nonumber \\
                     & = \frac{1}{\Gamma(\alpha)}\int_a^t  e^{-\lambda{t}} (t - u)^{\alpha-1} e^{\lambda{u}} \eta(u)du 
                       = \frac{1}{\Gamma(\alpha)}\int_a^t  e^{-\lambda(t-u)} (t - u)^{\alpha-1} \eta(u)du.
\label{eq:fracOU_0060}
\end{align}
Properties of FOU of Riemann-Liouville and Weyl types have been studied 
\cite{Eab06a,Eab06b,LimLITeo2007,LimLiTeo2008,LimTeo2009c}.
For the purpose of subsequent sections, only FOU of Weyl type would be considered. One has for $\alpha > 1/2$,
\begin{align}
  X_{\alpha,\lambda}(t) & = \frac{1}{\Gamma(\alpha)}\int_{-\infty}^t e^{-\lambda(t-u)} (t - u)^{\alpha-1} \eta(u)du \nonumber \\
                     & = \frac{1}{\Gamma(\alpha)}\int_{-\infty}^\infty e^{-\lambda(t-u)_{+}} (t - u)_{+}^{\alpha-1} \eta(u)du ,
\label{eq:fracOU_0070}
\end{align}
where $(x)_{+}^\mu = \bigl(\max(x,0)\bigr)^\mu$, $0^0 = 1$.
(\ref{eq:fracOU_0070}) is the moving average representation FOU of Weyl type. 
For $\alpha > 1/2$, $X_{\alpha,\lambda}(t)$ is a centred Gaussian stationary process with the following covariance and variance
\begin{align}
  C_{\alpha,\lambda}(t-s) & = \Bigl<X_{\alpha,\lambda}(t)X_{\alpha,\lambda}(s)\Bigr> 
                         =  \frac{1}{\sqrt{\pi}\Gamma(\alpha)}\left(\frac{|t-s|}{2\lambda}\right)^{\alpha-1/2} K_{\alpha-1/2}\bigl(\lambda|t-s|\bigr) , 
\label{eq:fracOU_0080}\\
  \sigma_{\alpha,\lambda}^2(t) & = \Bigl<\bigl(X_{\alpha,\lambda}(t)\bigr)^2\Bigr> 
                             = \frac{\Gamma(2\alpha-1)}{\bigl(\Gamma(\alpha)\bigr)^2(2\lambda)^{2\alpha-1}} .
\label{eq:fracOU_0090}
\end{align}

The spectral density of $X_{\alpha,\lambda}(t)$ is
\begin{align}
  S_{\alpha,\lambda}(k) & = \frac{1}{2\pi}\int_{\mathbb{R}}C_{\alpha,\lambda}(\tau)e^{ik\tau} d\tau = \frac{1}{2\pi}\frac{1}{\bigl(k^2 + \lambda^2\bigr)^\alpha}.
\label{eq:fracOU_0100}
\end{align}
FOU has the following spectral representation
\begin{align}
  X_{\alpha,\lambda}(t) & = \frac{1}{\sqrt{2\pi}}\int_{-\infty}^\infty \frac{e^{ikt}\widetilde{\eta}(k)dk}{(-ik+\lambda)^\alpha}, &
                     & \alpha > 1/2.
\label{eq:fracOU_0110}
\end{align}

Below are some properties of $X_{\alpha,\lambda}(t)$ which are relevant to the next section where TFBM will be treated as RFOU.

\noindent
\subsection{Scaling property}
\noindent
It is well-known that stationary Gaussian process cannot be a self-similar process 
\cite{SamorodnitskyTaqqu1994}.
$X_{\alpha,\lambda}(t)$ satisfies a weaker property of local self-similarity at small time scales 
\cite{KentWood1997,Adler1981,BenassiCohenIstas2003,AbramowitzStegun64}.
A Gaussian stationary process $Z(t)$ is locally self-similar of order $\kappa$ if its covariance $C(\tau)$  satisfies  
\begin{align}
  C(\tau) & = C(0) - A|\tau|^\kappa\bigl(1+o(1)\bigr) &
          & \text{as} \ |\tau| \to 0,
\label{eq:fracOU_0120}
\end{align}
where $0 < \kappa < 2$, and $A$ is a positive constant. 
One can show that FOU satisfies the local self-similarity by using the property of the modified Bessel function of second kind $K_\nu(z)$ 
for $|z| \to 0$ 
\cite{Falconer2003a}.
The small-time limit of the covariance of $X_{\alpha,\lambda}(t)$  behaves as
\begin{align}
  \frac{1}{\sqrt{\pi}\Gamma(\alpha)}\left(\frac{|\tau|}{2\lambda}\right)^{\alpha-1/2} K_{\alpha-1/2}\bigl(\lambda|\tau|\bigr)
     & \sim \frac{\Gamma(2\alpha-1)}{\bigl(\Gamma(\alpha)\bigr)^2(2\lambda)^{2\alpha-1}}
       + \frac{|\tau|^{2\alpha-1}}{2\Gamma(2\alpha)cos(\alpha\pi)}
       + o\bigl(\tau^2\bigr), &
    & |\tau| \to 0.
\label{eq:fracOU_0130}
\end{align}

Note that the class of Gaussian processes which satisfy 
(\ref{eq:fracOU_0120})
is also known as Adler processes 
\cite{Adler1981,BenassiCohenIstas2003},
which include 
the Gaussian stationary processes with stretched exponent covariance 
\cite{LimMunuandy2003}
and generalized Cauchy covariance 
\cite{LimTeo2009d}.
Instead of 
(\ref{eq:fracOU_0120}),
one can also adopt the definition of locally asymptotically self-similar property first introduced for multifractional Brownian motion 
\cite{BenassiJaffardRoux1997,BenassiCohenIstas2003,Falconer2003a}.
A stochastic process $Z(t)$ is locally asymptotically self-similarity at a point $t_\circ$ with order $\kappa$ if
\begin{align}
  \lim_{\epsilon \to 0}\biggl[\frac{Z(t_\circ + \epsilon{u}) - Z(t_\circ)}{\epsilon^\kappa}\biggr] & \ \hat{=} \ T_{t_\circ}(u),
\label{eq:fracOU_0140}
\end{align}
where $T_{t_\circ}(u)$ is a non-degenerate tangent process. 
Here the convergence is in the sense of finite dimensional distributions, and
$\hat{=}$
denotes equality in finite dimensional distributions.
Falconer 
\cite{Falconer2003a}
has shown that if the tangent process for a Gaussian process exists and is non-degenerate, 
then it is a self-similar Gaussian process with stationary increments. 
Since up to a multiplicative constant FBM is the only Gaussian self-similar process with stationary increments 
\cite{SamorodnitskyTaqqu1994},
the tangent process (\ref{eq:fracOU_0140}) is a FBM $B_\kappa(u)$ with Hurst index  $\kappa$.
One can easily verify (\ref{eq:fracOU_0140}) for $X_{\alpha,\lambda}(t)$ by direct computation using 
(\ref{eq:fracOU_0130}).
Thus, FOU behaves like FBM in small time scales.

In addition, $X_{\alpha,\lambda}(t)$ also satisfies the global scaling property similar to TFBM 
\cite{MeerschaertSbzikar2013,MeerschaertSabzikar2014}:
\begin{align}
  X_{\alpha,\lambda}(rt) & \ \hat{=} \ r^{\alpha-1/2} X_{\alpha,r\lambda}(t),
\label{eq:fracOU_0150}
\end{align}
where $r$ is a positive constant, and $X_{\alpha,r\lambda}(rt)$ is the same process as $X_{\alpha,\lambda}(t)$ with $\lambda$ replace by $r\lambda$.
(\ref{eq:fracOU_0150})
can be easily verified as follows. 
\begin{align}
  \Bigl<X_{\alpha,\lambda}(rt)X_{\alpha,\lambda}(rs)\Bigr> & = \frac{1}{\sqrt{\pi}\Gamma(\alpha)}\left(\frac{r|t-s|}{2\lambda}\right)^{\alpha-1/2} K_{\alpha-1/2}\bigl(r\lambda|t-s|\bigr) \nonumber \\
  & = \frac{r^{2\alpha-1}}{\sqrt{\pi}\Gamma(\alpha)}\left(\frac{|t-s|}{2r\lambda}\right)^{\alpha-1/2} K_{\alpha-1/2}\bigl(r\lambda|t-s|\bigr) \nonumber \\
  & = r^{2\alpha-1}\Bigl<X_{\alpha,r\lambda}(t)X_{\alpha,r\lambda}(s)\Bigr>.
\label{eq:fracOU_0160}
\end{align}

\subsection{Fractal or Hausdorff dimension}
\noindent
For the determination of the fractal dimension of the process $X_{\alpha,\lambda}(t)$, 
it is necessary to consider the local property of the process.  
First, we consider the H\"{o}lderian property of the sample path of $X_{\alpha,\lambda}(t)$.
A function $f:[a,b] \to \mathbb{R}$ is H\"{o}lderian of order $\kappa \in (0,1]$ if
\begin{align}
  \bigl|f(t) - f(s)\bigr| < K|t-s|^\kappa
\label{eq:fracOU_0170}
\end{align}
for all $s,t \in [a,b]$ for some constant $K > 0$.   
Since $\sigma_{\alpha,\Delta\tau}(t)$, the variance of the increment process of $X_{\alpha,\lambda}(t)$, satisfies
\begin{align}
  \sigma^2_{\alpha,\Delta\tau}(t) & = \Bigl<\bigl(X_{\alpha,\lambda}(t+\tau) - X_{\alpha,\lambda}(t)\bigr)^2\Bigr> \leq A|\tau|^{2\alpha-1},
\label{eq:fracOU_0180}
\end{align}
thus, almost surely the sample path of $X_{\alpha,\lambda}(t)$ is H\"{o}lderian of order $(\alpha-1/2)- \epsilon$ for all $\epsilon >0$.
For a process which is locally asymptotically self-similar of order $\kappa > 0$
and its sample paths are a.s. $\kappa -\epsilon$-H\"{o}lderian for all $\epsilon > 0$,  
then the Hausdorff dimension of its graph is a.s. equals to $2-\kappa$
\cite{KentWood1997,Adler1981,BenassiCohenIstas2003,Falconer2003b}.
Applying this result to $X_{\alpha,\lambda}(t)$ gives the Hausdorff dimension $5/2 - \alpha$
or $2 - H$ if $\alpha = H + 1/2$.
Thus, both the FOU and fractional Brownian motion have the same fractal dimension. 
This does not come as a surprise since FOU behaves like FBM locally.

\subsection{Short-range and long-range dependence}
\noindent
The usual way of characterising the memory of a stationary process in the time domain is in terms of decay rates of long-lag covariances, 
or in the frequency domain in terms of rates of divergence of spectral densities at low frequencies 
\cite{Beran1994}.
Note that for the covariance with large-time behavior
$C(\tau) \sim \tau^{-\beta}$, $\beta \in (0,1)$ as $\tau \to \infty$,
then the process is long range dependent (LRD). 
It can also be characterized by the power law divergence at 
the origin $S(k) \sim |k|^{\beta-1}$, $\beta \in (0,1)$ as $|k| \to 0$.
One can also use the following criterion for LRD property of a stationary Gaussian process 
\mbox{\cite{AyacheCohenLevyVehel2000,FandrinPAmblard2003}}.
The process is said to be LRD if its covariance $C(\tau)$ satisfies 
\begin{align}
  \int_0^\infty \bigl|C(\tau)\bigr| d\tau & = \infty .
\label{eq:fracOU_0190}
\end{align}
If the integral 
(\ref{eq:fracOU_0190})
is finite, the process is short range dependent (SRD). 

Note that for FOU, one has                                                         
\begin{align}
  \int_0^\infty \bigl|C(\tau)\bigr| d\tau & = \frac{1}{\sqrt{\pi}\Gamma(\alpha)} \int_0^\infty \left(\frac{|\tau|}{2\lambda}\right)^{\alpha-1/2}
                                            K_{\alpha-1/2}\bigl(\lambda|\tau|\bigr)d\tau
                                           = \frac{1}{2\lambda^{2\alpha}} .
\end{align}
Hence, FOU is a short memory process.

Consider a process defined by $Q(t) = \int_{-\infty}^\infty G(t-u)\eta(u)du$, where $G(t)$ is the response function. 
For the 3 cases  $G(t) = e^{-\lambda{t}},  e^{-\lambda{t}} t^{\alpha-1}, t^{\alpha-1}$, 
the resulting stochastic processes are respectively  
(i) Markov or memory-less Ornstein-Uhlenbeck process; 
(ii) short memory FOU; and 
(iii) long memory process FBM. 
Note that in case (ii) the exponential damping term in the response function prevents FOU to be LRD.


\section{Tempered Fractional Brownian Motion as Reduced Fractional Ornstein-Uhlenbeck Process}
\label{sec:TFBMRFOU}
\noindent
This section treats TFBM as reduced fractional Ornstein-Uhlenbeck process (RFOU). Recall that FBM defined by the Weyl integral
\begin{align}
  B_H(t) & = \frac{1}{\Gamma(H+1/2)}\int_{-\infty}^t (t-u)^{H-1/2} \eta(u)du
\label{eq:TFBMRFOU_0010}
\end{align}
is divergent. To ensure convergence, 
Mandelbrot and van Ness 
\cite{MandelbrotVanNess68}
modified 
(\ref{eq:TFBMRFOU_0010})
to a reduced FBM
\begin{align}
  B_H(t) - B_H(0) & = \frac{1}{\Gamma(H+1/2)}\Biggl[\int_{-\infty}^t (t-u)^{H-1/2} \eta(u)du - \int_{-\infty}^0 (-u)^{H-1/2} \eta(u)du\Biggr] ,
\label{eq:TFBMRFOU_0020}
\end{align}
which has since been generally accepted as the standard FBM. 

In contrast, FOU given by (\ref{eq:fracOU_0060}) is well-defined based on the Weyl-fractional integral. 
However, the covariance and variance of FOU of Weyl type is divergent in the ``massless'' limit $\lambda \to 0$.
As an illustration, one considers the simple case of ordinary Orsntein-Uhlenbeck process of Weyl type 
$X_{1,\lambda}(t)$ with covariance  
$\bigl<X_{1,\lambda}(t)X_{1,\lambda}(s)\bigr> = \bigl(e^{-\lambda|t - s|}\bigr)/2\lambda$
which diverges in the limit 
$\lambda \to 0$.
Its reduced process $X_{1,\lambda}(t) - X_{1,\lambda}(0)$ 
has the covariance
$\bigl(e^{-\lambda|t - s|} - e^{-\lambda|t|} - e^{-\lambda|s| +1}\bigr)/(2\lambda)$
with $\lambda \to 0$ limit given by
$\bigl(|t|+|s| - |t - s|\bigr)/2$ or $t\wedge s$,
which is just the covariance of Brownian motion. A similar situation exists for FOU, though it is more complicated 
(see 
\cite{LimTeo2007}
for details). 
The divergence of the covariance and variance in the $\lambda \to 0$ limit disappears if one considers the reduced process of FOU, 
with FBM as its limiting process.

The moving average definition of RFOU, denoted by $B_{\alpha,\lambda}(t)$, is defined for $\alpha > 1/2$:
\begin{align}
   B_{\alpha,\lambda}(t) & = X_{\alpha,\lambda}(t) - X_{\alpha,\lambda}(0)  \nonumber \\
                      & = \frac{1}{\Gamma(\alpha)} \int_{-\infty}^t e^{\lambda(t-u)} (t - u)^{\alpha-1} \eta(u) du
                       - \frac{1}{\Gamma(\alpha)} \int_{-\infty}^0 e^{\lambda(-u)} ( - u)^{\alpha-1} \eta(u) du .
\label{eq:TFBMRFOU_0030}
 \end{align}
From 
(\ref{eq:fracOU_0080})
and 
(\ref{eq:fracOU_0090}),
the covariance of $B_{\alpha,\lambda}(t)$ can be obtained directly as
\begin{align}
  \breve{C}_{\alpha,\lambda}(t,s) & = \Bigl<B_{\alpha,\lambda}(t)B_{\alpha,\lambda}(s)\Bigr> \nonumber \\
                               & = \Bigl<X_{\alpha,\lambda}(t)X_{\alpha,\lambda}(s)\Bigr> 
                                    - \Bigl<X_{\alpha,\lambda}(t)X_{\alpha,\lambda}(0)\Bigr>
                                    - \Bigl<X_{\alpha,\lambda}(s)X_{\alpha,\lambda}(0)\Bigr>
                                    + \Bigl<X_{\alpha,\lambda}(0)X_{\alpha,\lambda}(0)\Bigr> \nonumber \\
                               & = \frac{1}{\sqrt{\pi}\Gamma(\alpha)}\Biggl[
                                      \left(\frac{|t-s|}{2\lambda}\right)^{\alpha-1/2}
                                      K_{\alpha-1/2}(\lambda|t-s|)
                                    - \left(\frac{|t|}{2\lambda}\right)^{\alpha-1/2}
                                      K_{\alpha-1/2}(\lambda|t|)  \nonumber \\
                               & \qquad\qquad\qquad\qquad
                                    - \left(\frac{|s|}{2\lambda}\right)^{\alpha-1/2}
                                      K_{\alpha-1/2}(\lambda|s|)
                                    \Biggr]
                                 + \frac{\Gamma(2\alpha-1)}{\bigl(\Gamma(\alpha)\bigr)^2(2\lambda)^{2\alpha-1}}
\label{eq:TFBMRFOU_0040}
\end{align}
By letting $H = \alpha -1/2$, and using
\begin{align}
  c_t & = \frac{2\Gamma(2H)}{\bigl(\Gamma(H+1/2)\bigr)^2(2\lambda)^{2H}} - \frac{2}{\sqrt{\pi}\Gamma(H+1/2)}\left(\frac{1}{2\lambda|t|}\right)^H K_H(\lambda|t|),
\label{eq:TFBMRFOU_0050}
\end{align}
the covariance of the RFOU becomes 
\begin{align}
  \breve{C}_{\alpha,\lambda}(t,s) & = \frac{1}{2}\Bigl[c_t|t|^{2H} + c_s|s|^{2H} - c_{t-s}|t-s|^{2H}\Bigr],
\label{eq:TFBMRFOU_0060}
\end{align}
which is just the covariance of TFBM [1,2] up to a multiplicative constant $\bigl(\Gamma(H+1/2)\bigr)^{-2}$
due to the additional $1/\Gamma(\alpha)$ term in the definition of $B_{\alpha,\lambda}(t)$.
The variance of $B_{\alpha,\lambda}(t)$ is
\begin{align}
  \breve{\sigma}_{\alpha,\lambda}^2(t) & = \Bigl<\bigl(B_{\alpha,\lambda}(t)\bigr)^2\Bigr>
                                     =  \frac{2\Gamma(2H)}{\bigl(\Gamma(H+1/2)\bigr)^2(2\lambda)^{2H}}   
                                     -  2\frac{1}{\sqrt{\pi}\Gamma(H+1/2)}\left(\frac{|t|}{2\lambda}\right)^H K_H(\lambda|t|).
\label{eq:TFBMRFOU_0070}
\end{align}

\begin{figure}[ht]
  \centering
  \fbox{
  \includegraphics[bb = 0 0 350 280,width = 0.5\textwidth]{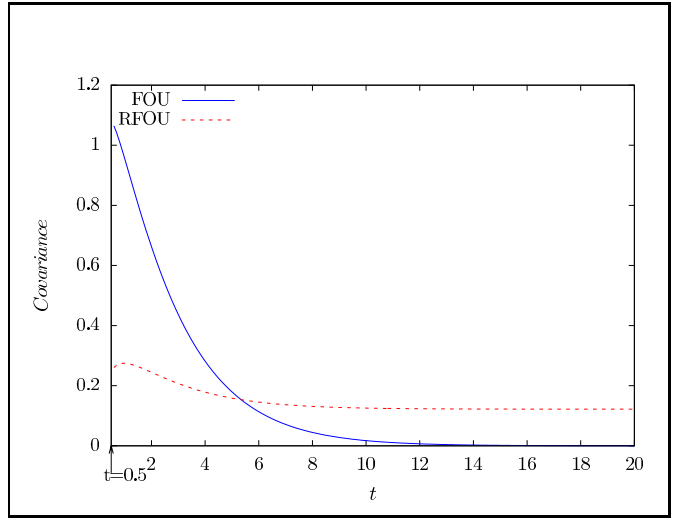}
  }
  \caption{covaraince of FOU and RFOU: $s=0.5$, $\lambda = 0.5, H = 0.75$}
  \label{fig:FOUvsTFBM_0010}
\end{figure}

Figure \mbox{\ref{fig:FOUvsTFBM_0010}} shows an example of the covariance functions of FOU and RFOU. 
Note that 
(\ref{eq:TFBMRFOU_0070})
can be expressed as $\breve{\sigma}^2_{\alpha,\lambda}(t) = 2\bigl(\sigma^2_{\alpha,\lambda}(t) - C_{\alpha,\lambda}(t)\bigr)$ and
$C_{\alpha,\lambda}(t)$ decays from $\sigma^2_{\alpha,\lambda}(t)$ to zero (see \mbox{(\ref{eq:fracOU_0080})} and graph in Figure \mbox{\ref{fig:FOUvsTFBM_0010}})
one obtains $\sigma^2_{\alpha,\lambda}(t) \leq \breve{\sigma}^2_{\alpha,\lambda}(t) \leq 2\sigma^2_{\alpha,\lambda}(t)$.

Many of the basic the properties of RFOU or TFBM are direct consequences or some modifications of the properties of FOU. 
These include the following.

\subsection{Spectral representation}
\noindent
Using the spectral representation of $X_{\alpha,\lambda}(t)$ given by 
(\ref{eq:fracOU_0110})
one gets the following harmonizable representation for TFBM:
\begin{align}
  B_{\alpha,\lambda}(t) & 
                       = \frac{1}{\sqrt{2\pi}} \int_{-\infty}^\infty \frac{\Bigl(e^{-kt} - 1\Bigr)\widetilde{\eta}(k)dk}{(-ik + \lambda)^\alpha}.
\label{eq:TFBMRFOU_0080}
\end{align}
The covariance of $B_{\alpha,\lambda}(t)$ has the following harmonizable representation
\begin{align}
  \breve{C}_{\alpha,\lambda}(t,s) & = \frac{1}{2\pi} \int_{-\infty}^\infty \frac{\Bigl(e^{ikt} - 1\Bigr)\Bigl(e^{-iks} - 1\Bigr)}{\bigl(|k|^2 + \lambda^2\bigr)^\alpha} dk  
                                = \frac{1}{2\pi} \int_{-\infty}^\infty \frac{1+e^{ik(t-s)} - e^{ikt} - e^{-iks}}{\bigl(|k|^2 + \lambda^2\bigr)^\alpha} dk,
\label{eq:TFBMRFOU_0090}
\end{align} 
which can be evaluated to give (\ref{eq:TFBMRFOU_0020}).

\subsection{Stationary increments}
\noindent
In contrast to FOU, reduced process of FOU or TFBM is non-stationary. 
However, TFBM is asymptotically stationary, 
that is $B_{\alpha,\lambda}(t)$  is stationary at very large $t$.
Consider \mbox{(\ref{eq:TFBMRFOU_0040})} the covariance of TFBM $\breve{C}_{\alpha,\lambda}(t,s)$ for $t, s \to \infty$,
and $|t-s|$ finite.
Since $K_\nu(z) \sim \sqrt{\frac{\pi}{2z}}e^{-z}$ as $z \to \infty$, one has
\begin{align}
   \breve{C}_{\alpha,\lambda}(t,s)  &  \sim  \left(\frac{|t-s|}{2\lambda}\right)^{\alpha-1/2} K_{\alpha - 1/2}\bigl(\lambda|t-s|\bigr)
                                        + \frac{\Gamma(2\alpha)}{\bigl(\Gamma(\alpha)\bigr)^2(2\lambda)^{2\alpha -1}}.
\label{eq:TFBMRFOU_0100}
\end{align}
On the other hand, the increment process of TFBM, $\Delta_\tau B_{\alpha,\lambda}(t) = B_{\alpha,\lambda}(t+\tau) - B_{\alpha,\lambda}(t)$, 
is stationary. 
One has the tempered fractional Gaussian noise
\begin{align}
  \Delta_\tau B_{\alpha,\lambda}(t) & = \bigl(X_{\alpha,\lambda}(t+\tau) - X_{\alpha,\lambda}(0)\bigr) - \bigl(X_{\alpha,\lambda}(t) - X_{\alpha,\lambda}(0)\bigr)   \nonumber \\
       & = X_{\alpha,\lambda}(t+\tau) - X_{\alpha,\lambda}(t) = \Delta_\tau X_{\alpha,\lambda}(t).
\label{eq:TFBMRFOU_0110}
\end{align}
Since $X_{\alpha,\lambda}(t)$ is a stationary process, so is its increment process $\Delta_\tau X_{\alpha,\lambda}(t)$, hence $\Delta_\tau B_{\alpha,\lambda}(t)$.

From (\ref{eq:TFBMRFOU_0080}),
the spectral representation of the increment process is
\begin{align}
 \Delta_\tau X_{\alpha,\lambda}(t) & = \frac{1}{\sqrt{2\pi}} \int_{-\infty}^\infty 
                                       \frac{\Bigl(e^{-ik(t+\tau)} - e^{-ikt}\Bigr)\widetilde{\eta}(k)dk}
                                       {(-ik + \lambda)^\alpha} ,
\label{eq:TFBMRFOU_0120}
\end{align}
which gives the covariance of the increment process
\begin{align}
  \breve{C}^\tau_{\alpha,\lambda}(t,s) & = \bigl<\Delta_\tau B_{\alpha,\lambda}(t)\Delta_\tau B_{\alpha,\lambda}(s)\bigr> \nonumber \\
& = \frac{1}{2\pi} \int_{-\infty}^\infty \frac{\Bigl(e^{-ik(t+\tau)} - e^{-ikt}\Bigr)\Bigl(e^{ik(s+\tau)} - e^{iks}\Bigr)dk}{\bigl(|k|^2 + \lambda^2\bigr)^\alpha} \nonumber \\
   & = \frac{1}{\sqrt{\pi}\Gamma(\alpha)} \Biggl[
        2\left(\frac{|t-s|}{2\lambda}\right)^{\alpha - 1/2} K_{\alpha-1/2}\bigl(\lambda|t-s|\bigr)  \nonumber \\
   & \qquad\qquad\qquad
       - \left(\frac{|t-s+\tau|}{2\lambda}\right)^{\alpha - 1/2} K_{\alpha-1/2}\bigl(\lambda|t-s+\tau|\bigr) \nonumber \\
   & \qquad\qquad\qquad
       - \left(\frac{|t-s-\tau|}{2\lambda}\right)^{\alpha - 1/2} K_{\alpha-1/2}\bigl(\lambda|t-s-\tau|\bigr)
     \Biggr] .
\label{eq:TFBMRFOU_0130}
\end{align}

Another way to express the covariance of the increment process of TFBM is as follows.
\begin{align}
  \breve{C}_{\alpha,\lambda}^\tau (t,s) & = \frac{1}{2\pi} \int_{-\infty}^\infty \frac{2-2\cos(k\tau)}{\bigl(|k|^2 + \lambda^2\bigr)^\alpha}e^{ik(t-s)}dk \nonumber \\
  & = \frac{2}{\pi} \int_{-\infty}^\infty \frac{\sin^2(k\tau)}{\bigl(|k|^2 + \lambda^2\bigr)^\alpha}e^{ik(t-s)}dk,
\end{align}
which gives the spectral density of $\Delta_\tau B_{\alpha,\lambda}(t)$ as
$\breve{S}_{\alpha,\lambda}^\tau(k) = 4 \frac{\sin^2(k\tau)}{\bigl(|k|^2 + \lambda^2\bigr)^\alpha}$.
\\

Note that instead of the increment process (32), 
it is possible to define the fractional tempered Gaussian noise correspond to TFBM as the derivative of $B_{\alpha,\lambda}(t)$ in the sense of generalized functions, just like fractional Gaussian noise in terms of generalized derivative of FBM. 
One has  $\xi_{\alpha,\lambda}(t) = \frac{dB_{\alpha,\lambda}(t)}{dt}$ , which can be shown to be a stationary SRD process. 
See  \ref{sec:TFGN} for more details.

\subsection{Scaling property}
\noindent
FOU $X_{\alpha,\lambda}(t)$ satisfies the scaling property 
(\ref{eq:fracOU_0150})
for all $t$.
By applying this property to the covariance of RFOU 
(\ref{eq:TFBMRFOU_0040})
one can easily show that 
(\ref{eq:fracOU_0150})
also holds for $B_{\alpha,\lambda}(t)$ with
\begin{align}
  \Bigl<B_{H,\lambda}(rt)B_{H,\lambda}(rs)\Bigr> & = r^{2H}\Bigl<B_{H,\lambda}(rt)B_{H,\lambda}(rs)\Bigr> .
\label{eq:TFBMRFOU_0130}
\end{align}

The above scaling property can loosely be regarded as a generalization of self-similarity property. 
In additional to FOU, such a scaling property also holds for stationary process with appropriately normalised stretched exponential covariance 
$C_{\nu,\lambda}(t) = e^{-(\lambda|t|)^\nu}/(2\lambda^\nu)$,
$\nu \in (0,2]$, $\lambda > 0$.
This covariance diverges in the limit $\lambda \to 0$.
It is interesting to note that its reduced process also satisfies the scaling property, 
and becomes FBM as $\lambda \to 0$.
Thus, one may say that the scaling property 
(\ref{eq:fracOU_0150}) 
is not unique to the processes of Ornstein-Uhlenbeck type, 
which include FOU and TFBM, 
it is satisfied by a wider class of Gaussian stationary and non-stationary processes,
in particular Gaussian stationary processes of Adler class 
\cite{Adler1981,BenassiCohenIstas2003}
and their reduced processes.

Just like FOU, TFBM also satisfies the locally self-similarity with the condition 
\mbox{(\ref{eq:fracOU_0120})}
modified for non-stationary Gaussian process 
\cite{LimLITeo2007,KentWood1997}:
\begin{align}
  \breve{C}_{\alpha,\lambda}(t,t+\tau) & = \frac{1}{2} \bigl(\breve{C}_{\alpha,\lambda}(t,t) + \breve{C}_{\alpha,\lambda}(t+\tau,t+\tau)\bigr)
                                      - A|\tau|^\kappa\bigl(1+o(1)\bigr), &
                                    & 0 < \kappa < 2, \text{as} \ |\tau| \to 0.   
\label{eq:TFBMRFOU_0140}
\end{align}
Direct computation of $\wideparen{C}_\alpha(t,t+\tau)$  and using 
(\ref{eq:fracOU_0130})
one gets for the limit $|\tau| \to 0$,
\begin{align}
   \breve{C}_{\alpha,\lambda}(t,t+\tau) & = \frac{1}{2\Gamma(2\alpha)\cos(\alpha\pi)}
                                          \Bigl(
                                          |\tau|^{2\alpha-1} - |t|^{2\alpha-1} - |t+\tau|^{2\alpha-1} + o\bigl(\tau^{2\alpha-1}\bigr)
                                          \Bigr).
\label{eq:TFBMRFOU_0150}
\end{align}
Using $-\bigl(2\Gamma(2\alpha)\cos(\alpha\pi)\bigr)^{-1} = \Gamma(1-2H)\cos(H\pi)/(2H\pi)$, 
where $\alpha = H + 1/2$, one gets
\begin{align}
  \breve{C}_{\alpha,\lambda}(t,t+\tau) & = \frac{\Gamma(1-2H)\cos(H\pi)}{2H\pi}\Bigl[|t+\tau|^{2H} + |t|^{2H} - |\tau|^{2H} + o\bigl(\tau^{2H}\bigr)\Bigr].
\label{eq:TFBMRFOU_0160}
\end{align}
Similarly, one has in the limit  $|\tau| \to 0$,
\begin{align}
    \frac{1}{2}\Bigl(\breve{C}_{\alpha,\lambda}(t+\tau,t+\tau) - \breve{C}_{\alpha,\lambda}(t,t)\Bigr) & = \frac{\Gamma(1-2H)\cos(H\pi)}{2H\pi}\Bigl[|t+\tau|^{2H} + |t|^{2H}  + o\bigl(\tau^{2H}\bigr)\Bigr].
\label{eq:TFBMRFOU_0170}
\end{align}
So 
(\ref{eq:TFBMRFOU_0140})
holds and $B_{\alpha,\lambda}(t)$ is locally self-similar. 

The tangent process 
(\ref{eq:fracOU_0140})
can be determined by noting that the increment process of FOU is the same as the increment process of RFOU. 
Thus, both FOU and RFOU (or TFBM) behave locally like FBM. 
Note that this property is not specific to these two processes. 
Gaussian processes belonging to Adler class satisfy the same property
\cite{Adler1981,BenassiCohenIstas2003}.

\subsection{Fractal dimension}
\noindent
Fractal dimension of TFBM can be obtained in the same way as for FOU. 
In particular, 
(\ref{eq:fracOU_0180})
holds for RFOU so that the sample path of $B_{H,\lambda}(t)$ is H\"{o}lderian of order 
$(\alpha - 1/2) -\epsilon$ for all $\epsilon > 0$.
Hence the Hausdorff dimension of its graph is a.s. equals to $5/2 - \alpha$ or $2 - 2H$ for 
$\alpha = H + \nicefrac{1}{2}$, which is the same as the fractal dimension of FBM. 
Thus, the three processes, FOU, TFBM and FBM have the same fractal dimension as they have the same local behaviour.

\subsection{Long-range dependence}
\noindent
Unlike FOU which is a short memory process, TFBM is LRD.
Condition 
(\ref{eq:fracOU_0190})
can be generalised to the non-stationary Gaussian process 
\cite{AyacheCohenLevyVehel2000,FandrinPAmblard2003}.
In this case, instead of the covariance one considers the correlation function
\begin{align}
  \breve{R}_{\alpha,\lambda}(t,t+\tau) & = \frac{\breve{C}_{\alpha,\lambda}(t,t+\tau)}{\bigl[\breve{C}_{\alpha,\lambda}(t,t)\breve{C}_{\alpha,\lambda}(t+\tau,t+\tau)\bigr]^{1/2}} 
                                      = \frac{\breve{C}_{\alpha,\lambda}(t,t+\tau)}{\bigl[\breve{\sigma}_{\alpha,\lambda}^2(t)\breve{\sigma}_{\alpha,\lambda}^2(t+\tau)\bigr]^{1/2}} .
\label{eq:TFBMRFOU_0180}
\end{align}
A non-stationary Gaussian process is said to be LRD
if its correlation function $\breve{R}_{\alpha,\lambda}(t,t+\tau)$  satisfies
\begin{align}
  \int_0^\infty \bigl|\breve{R}_{\alpha,\lambda}(t,t+\tau\bigr|d\tau & = \infty,
\label{eq:TFBMRFOU_0190}
\end{align}
otherwise it is SRD.
Alternatively, the process is said to  be LRD if 
$\breve{R}_{\alpha,\lambda}(t,t+\tau) \sim \tau^\gamma$, 
$-1 < \gamma < 0$ for $\tau \to \infty$, for all $t > 0$.

From 
(\ref{eq:TFBMRFOU_0100})
one has for $\tau \to \infty$, $\breve{C}_{\alpha,\lambda}(t+\tau,t) \sim C_{\alpha,\lambda}(0) - C_{\alpha,\lambda}(t) = \breve{\sigma}_{\alpha,\lambda}^2(t)/2$.
From 
(\ref{eq:fracOU_0090})
and 
(\ref{eq:TFBMRFOU_0070}),
and 
$K_\nu(\tau) \sim \sqrt{\pi/2\tau}e^{-\tau} \to 0$ as $\tau \to \infty$, 
one gets $\breve{\sigma}_{\alpha,\lambda}^2(t+\tau) \to \sigma_{\alpha,\lambda}^2(t)$.  
Thus, one has for $\tau \to \infty$,
\begin{align}
  \breve{R}_{\alpha,\lambda}(t,t+\tau) 
                                      \sim  \frac{1}{2} \sqrt{\frac{\breve{\sigma}_{\alpha,\lambda}^2(t)}{\sigma_{\alpha,\lambda}^2(t)}}   
                                      > 0,   
\label{eq:TFBMRFOU_0200}
\end{align}
which implies the TFBM is a long memory process.


\section{Mixed Tempered Fractional Brownian Motion}
\label{sec:MixTFBM}
\noindent
Mixed TFBM can be defined in a similar way as mixed FBM
\cite{Cheridito2001,El-Nouty2003}.
Denote mixed TFBM by $B_{\boldsymbol{\alpha},\boldsymbol{\lambda}}(t)$ with
$\boldsymbol{\alpha} = \bigl(\alpha_i, i = 1, \cdots, n\bigr)$,
$\boldsymbol{\lambda} = \bigl(\lambda_i, i = 1, \cdots, n\bigr)$,
\begin{align}
  B_{\boldsymbol{\alpha},\boldsymbol{\lambda}}(t) & = \sum_{i=1}^n b_i B_{\alpha_i,\lambda_i}(t), &
                                b_i & > 0, 
\label{eq:MixTFBM_0010}
\end{align}
and $B_{\alpha_i,\lambda_i}(t)$, $i = 1, \cdots, n$, $t \in [0,\infty)$ are $n$ independent TBFM 
with distinct indices $\alpha_i$ and constant mixture coefficients $b_i$.

$B_{\boldsymbol{\alpha},\boldsymbol{\lambda}}(t)$ is a centred Gaussian process with covariance given by
\begin{align}
  \breve{C}_{\boldsymbol{\alpha},\boldsymbol{\lambda}}(t,s) & = \Bigl<B_{\boldsymbol{\alpha},\boldsymbol{\lambda}}(t)B_{\boldsymbol{\alpha},\boldsymbol{\lambda}}(s)\Bigr>
                                              =  \sum_{i,j=1}^n b_ib_j\Bigl<B_{\alpha_i,\lambda_i}(t)B_{\alpha_i,\lambda_i}(s)\Bigr>  \nonumber \\
       & = \frac{1}{2} \sum_{i=1}^n b_i^2 \Bigl[
            c_{i,t}|t|^{2H_i}
           + c_{i,s}|s|^{2H_i}
           - c_{i,t-s} |t-s|^{2H_i}
         \Bigr] \nonumber \\
       & = \sum_{i=1}^n b_i^2\breve{C}_{\alpha_i,\lambda_i}(t,s) ,
\label{eq:MixTFBM_0020}
\end{align}
where
\begin{align}
  c_{i,t} & = \frac{\Gamma(2H_i)}{\bigl(\Gamma(H_i+1/2)\bigr)^2(2\lambda_i)^{2H_i}}
                       - \frac{2}{\sqrt{\pi}\Gamma(H_i+1/2)}\left(\frac{1}{2\lambda_i|t|}\right)^{H_i} K_{H_i}\bigl(\lambda_i|t|\bigr),
\label{eq:MixTFBM_0030}
\end{align}
and $\alpha_i = H_i + 1/2$.

Its variance is  
\begin{align}
  \Bigl<\bigl(B_{\boldsymbol{\alpha},\boldsymbol{\lambda}}(t)\bigr)^2\Bigr> & = \sum_{i=1}^n b_i^2 \Bigl<\bigl(B_{\alpha_i,\lambda_i}(t)\bigr)^2\Bigr> \nonumber \\
     & =  \sum_{i=1}^n 2b_i^2 
          \Biggl(\frac{\Gamma(2H_i)}{\bigl(\Gamma(H_i+1/2)\bigr)^2(2\lambda_i)^{2H_i}}
                       - \frac{1}{\sqrt{\pi}\Gamma(H_i+1/2)}\left(\frac{|t|}{2\lambda_i}\right)^{H_i} K_{H_i}\bigl(\lambda_i|t|\bigr)
          \Biggr).
\label{eq:MixTFBM_0040}
\end{align}

Since mixed TFBM is a linear combination of $n$ independent copies of TFBM, its properties are direct consequences of the properties of TFBM.

\subsection{Mixed scaling}
\noindent
$B_{\boldsymbol{\alpha},\boldsymbol{\lambda}}(t)$
satisfies the following mixed scaling property:                              
\begin{align}
  B_{\boldsymbol{\alpha},\boldsymbol{\lambda}}(rt) & = \sum_{i=1}^n b_ir^{\alpha_i-1/2} B_{\alpha_i,r\lambda_i}(t).
\label{eq:MixTFBM_0050}
\end{align}
By noting that the scaling property 
(\ref{eq:TFBMRFOU_0130})
holds for each $B_{\alpha_i,\lambda_i}(t)$, $i = 1, \cdots, n$,
hence the linear combination of $n$ independent copies of $B_{\alpha_i,\lambda_i}(t)$ satisfies the mixed scaling property 
(\ref{eq:MixTFBM_0050}).
Since $B_{\boldsymbol{\alpha},\boldsymbol{\lambda}}(t)$ is a centred Gaussian process one only needs to prove that 
$\sum_{i=1}^n b_iB_{\alpha_i,\lambda_i}(rt)$ and $\sum_{i=1}^nr^{\alpha_i} b_iB_{\alpha_i,\lambda_i}(t)$
have the same covariance function. 

 $B_{\boldsymbol{\alpha},\boldsymbol{\lambda}}(t)$ satisfies the following mixed locally self-similarity property:
\begin{align}
  \breve{C}_{\boldsymbol{\alpha},\boldsymbol{\lambda}}(t,t+\tau) & = \frac{1}{2} \sum_{i=1}^n 
                                                          \Bigl(
                                                          \breve{C}_{\alpha_i,\lambda_i}(t,t)
                                                          + \breve{C}_{\alpha_i,\lambda_i}(t+\tau,t+\tau)
                                                          - A_i|\tau|^{2\alpha_i-1} \bigl(1+o(1)\bigr)
                                                          \Bigr), \nonumber \\
                                                  &  \hspace{0.65\textwidth} |\tau|  \to 0.                 
\label{eq:MixTFBM_0060}
\end{align}
This follows from the fact that each TFBM $B_{\alpha_i,\lambda_i}(t)$ independently satisfies locally self-similar property 
(\ref{eq:TFBMRFOU_0140}).

Consider the tangent process of $B_{\boldsymbol{\alpha},\boldsymbol{\lambda}}(t)$ at a point $t_\circ$:
\begin{align}
  \lim_{\epsilon\to 0} \Biggl[
                   \frac{B_{\boldsymbol{\alpha},\boldsymbol{\lambda}}(t_\circ+\epsilon{u}) - B_{\boldsymbol{\alpha},\boldsymbol{\lambda}}(t_\circ)}{\epsilon^{\boldsymbol{\kappa}}}
                     \Biggr]
                 & =   \lim_{\epsilon\to 0} \Biggl[
                       \sum_{i=1}^n b_i
                       \frac{B_{\alpha_i,\lambda_i}(t_\circ+\epsilon{u}) - B_{\alpha_i,\lambda_i}(t_\circ)}{\epsilon^{\kappa_i}}
                     \Biggr] .
\label{eq:MixTFBM_0070}
\end{align}
Since the tangent process of each $B_{\alpha_i,\lambda_i}(t)$ at a point $t_\circ$ is FBM $B_{H_i}(t_\circ)$ indexed by $H_i = \alpha_i - 1/2$,
it is straight forward to verify that  the tangent process of mixed TFBM 
(\ref{eq:MixTFBM_0070})
at a point $t_\circ$ is given by the mixed FBM $\sum_{i=1}^n b_i B_{H_i}(t_\circ)$.

\subsection{Stationary increments}
\noindent
Since each $B_{\alpha_i,\lambda_i}(t)$ has stationary increments, due to the independence of $B_{\alpha_i,\lambda_i}(t)$, $i = 1,\cdots, n$,
hence the increments of their linear combination are stationary. One has
\begin{gather}
  \Bigl<\bigl(B_{\boldsymbol{\alpha},\boldsymbol{\lambda}}(t+\tau)\bigr) - B_{\boldsymbol{\alpha},\boldsymbol{\lambda}}(t) \bigr) 
     \bigl(B_{\boldsymbol{\alpha},\boldsymbol{\lambda}}(s+\tau)\bigr) - B_{\boldsymbol{\alpha},\boldsymbol{\lambda}}(s) \bigr)
  \Bigr> \qquad\qquad\qquad\qquad\qquad\qquad \nonumber \\
  = \sum_{i,j=1}^n b_ib_j
  \Bigl<\bigl(B_{\alpha_i,\lambda_i}(t+\tau)\bigr) - B_{\alpha_i,\lambda_i}(t) \bigr) 
      \bigl(B_{\alpha_i,\lambda_i}(s+\tau)\bigr) - B_{\alpha_i,\lambda_i}(s) \bigr) 
  \Bigr>  \nonumber \\
   = \sum_{i,j=1}^n b_i^2
  \Bigl<\bigl(X_{H_i,\lambda_i}(t+\tau)\bigr) - X_{H_i,\lambda_i}(t) \bigr) 
      \bigl(X_{H_i,\lambda_i}(s+\tau)\bigr) - X_{H_i,\lambda_i}(s) \bigr) 
  \Bigr> \nonumber \\ 
   = \sum_{i,j=1}^n b_i^2
     \frac{1}{\sqrt{\pi}\Gamma(\alpha_i)} 
     \Biggl[
     2\left(\frac{|t-s|}{2\lambda_i}\right)^{\alpha_i-1/2} K_{\alpha_i-1/2}\bigl(\lambda_i|t-s|\bigr)\qquad\qquad\qquad \nonumber \\
     - \left(\frac{|t-s+\tau|}{2\lambda_i}\right)^{\alpha_i-1/2} K_{\alpha_i-1/2}\bigl(\lambda_i|t-s+\tau|\bigr) \qquad\qquad\nonumber \\
     - \left(\frac{|t-s-\tau|}{2\lambda_i}\right)^{\alpha_i-1/2} K_{\alpha_i-1/2}\bigl(\lambda_i|t-s-\tau|\bigr)
     \Biggr].
\label{eq:MixTFBM_0080}
\end{gather}

\subsection{Long range dependence}
\noindent
The sum of $n$ independent LRD Gaussian processes 
$B_{\alpha_i,r\lambda_i}(t)$, $i = 1, \cdots, n$,
is LRD. 
One can also verify this property by considering the correlation of $B_{\boldsymbol{\alpha},\boldsymbol{\lambda}}(t)$ is
\begin{align}
  \breve{R}_{\boldsymbol{\alpha},\boldsymbol{\lambda}}(t+\tau,t) & = \frac{\breve{C}_{\boldsymbol{\alpha},\boldsymbol{\lambda}}(t+\tau,t)}
                                                          {\sqrt{\Bigl<\bigl(B_{\boldsymbol{\alpha},\boldsymbol{\lambda}}(t+\tau)\bigr)^2\Bigr>
                                                                 \Bigl<\bigl(B_{\boldsymbol{\alpha},\boldsymbol{\lambda}}(t)\bigr)^2\Bigr>}} 
     = \frac{\sum_{i=1}^n \breve{C}_{H_i,\lambda_i}(t+\tau,t)}{\sqrt{\vphantom{\Bigl[}\sum_{i=1}^n\breve{\sigma}_i^2(t)\sum_{i=1}^n\breve{\sigma}_i^2(t+\tau)}} .
\label{eq:MixTFBM_0090}
\end{align}
$\breve{C}_{H_i,\lambda_i}(t+\tau,t) \to \frac{\breve{\sigma}_i^2(t)}{2}$ as $\tau \to \infty$.
Let $\breve{\sigma}_{\text{max}}^2 =  \max\{\breve{\sigma}_i^2 , i = 1, \cdots, n\}$ and
 $\breve{\sigma}_{\text{min}}^2 =  \min\{\breve{\sigma}_i^2 , i = 1, \cdots, n\}$.
Therefore,
\begin{align}
  \breve{R}_{\boldsymbol{\alpha},\boldsymbol{\lambda}}(t+\tau,t & > \frac{\breve{\sigma}_{\text{min}}^2(t)}{2\sqrt{\breve{\sigma}_{\text{max}}^2(t+\tau)\breve{\sigma}_{\text{max}}^2(t)}}
\sim
\frac{\breve{\sigma}_{\text{min}}^2(t)}{2\breve{\sigma}_{\text{max}}^2(t)}
  > 0 .
\label{eq:MixTFBM_0100}
\end{align}
which implies $B_{\boldsymbol{\alpha},\boldsymbol{\lambda}}(t)$ is long-range dependent.

\subsection{Fractal dimension}
\noindent
The Hausdorff dimension of the graph of $B_{\alpha_i,\lambda_i}(t)$ is $5/2 - \alpha_i$, with probability 1. 
Based on a standard results \cite{Falconer2003b} the fractal dimension of $B_{\boldsymbol{\alpha},\boldsymbol{\lambda}}(t)$ as
$5/2- \min\limits_{1 \leq i \leq n} \alpha_i$

Mixed TFBM can be used to model systems which require variable tempering factor $\lambda_i$, or different index $\alpha_i$, 
or both at different time intervals.
For example, for small time scales $B_{\boldsymbol{\alpha},\boldsymbol{\lambda}}(t)$  behaves like mixed FBM, $B_{\alpha_i,\lambda_i}(t)$ 
with smaller ``Hurst index'' dominates. In the intermediate time scales, there would be an interplay between the tempering factor $\lambda_i$
and index $\alpha_i$.
Finally, for large time scales, it is dictated by $B_{\alpha_i,\lambda_i}(t)$ with higher tempering factor $\lambda_i$ which enters as an exponential term.


\section{Tempered Fractional Brownian Motion with Two Indices}
\label{sec:2indRFOU}
\noindent
For TFBM considered in section 
\ref{sec:TFBMRFOU},
the long time and short time properties of the process are characterized by the same index $\alpha$. 
It will be useful from practical point of view that TFBM can be generalised to a process 
such that its long and short time behavior can be described by two separate parameters. 
One possible way to achieve this goal is to extend FOU to two indices 
\cite{LimLiTeo2008,LimTeo2009c} 
and consider its reduced process.  

The fractional Langevin equation 
\mbox{(\ref{eq:fracOU_0050})}
can be extended to two indices 
\begin{align}
  \bigl({_{-\infty}D_t^\beta} + \lambda^\beta\bigr)^\alpha X_{\alpha\beta,\lambda}(t) & = \eta(t), &
                                                                 \alpha\beta & > 1/2 .
\label{eq:2indRFOU_0010}
\end{align}
Here $\lambda$ is replaced by $\lambda^\beta$ to preserve the scaling property. 
The shifted fractional differential operator with two indices $\bigl({_{-\infty}D_t^\beta} + \lambda^\beta\bigr)^\alpha$ can be expanded as binomial series:
\begin{align}
  \bigl({_{-\infty}D_t^\beta} + \lambda^\beta\bigr)^\alpha & = \sum_{j=0}^\infty \binom{\alpha}{j} \lambda^{\beta{j}} {_{-\infty}D_t^{\beta(\alpha - j)}}.
\label{eq:2indRFOU_0020}
\end{align}
There does not exist a nice operator identity similar to
(\ref{eq:fracOU_0040})
as in case of the shifted fractional derivative of single index.

By using Fourier transform method, the solution of 
(\ref{eq:2indRFOU_0010})
is found to be a stationary Gaussian process
\begin{align}
   X_{\alpha\beta,\lambda}(t) & = \frac{1}{\sqrt{2\pi}}\int_{-\infty}^\infty \frac{e^{ikt}\widetilde{\eta}(k)}{\bigl((ik)^\beta + \lambda^\beta\bigr)^\alpha} dk.
\label{eq:2indRFOU_0030}
\end{align}
The condition $\alpha\beta > 1/2$ is necessary to ensure the above integral is finite. 
This process has a more complicated spectral density,
\begin{align}
  S_{\alpha\beta}(k) & = \frac{1}{2\pi\bigl|(ik)^\beta + \lambda^\beta\bigr|^{2\alpha}}
                    = \frac{1}{2\pi\bigl(|k|^{2\beta} +  2\lambda^\beta|k|^\beta\cos(\alpha\pi/2)+\lambda^{2\beta}\bigr)^\alpha} .  
\label{eq:2indRFOU_0040}
\end{align}
The inverse Fourier transform of 
(\ref{eq:2indRFOU_0040})
in general cannot be evaluated to give a closed analytic expression of the covariance for $X_{\alpha\beta,\lambda}(t)$.
Despite this, many of the basic properties of $ X_{\alpha\beta,\lambda}(t)$  can still be obtained and studied 
\cite{LimLiTeo2008,LimTeo2009c}.
Hence, one can consider the reduced process associated with $X_{\alpha\beta,\lambda}(t)$ given by
\begin{align}
   B_{\alpha\beta,\lambda}(t) & = X_{\alpha\beta,\lambda}(t) - X_{\alpha\beta,\lambda}(0)
                            = \frac{1}{\sqrt{2\pi}}\int_{-\infty}^{\infty} \frac{\bigl(e^{ikt} - 1\bigr)\widetilde{\eta}(k)}{\bigl((ik)^\beta + \lambda^\beta\bigr)^\alpha} dk,
\label{eq:2indRFOU_0050}
\end{align}
and examine its properties just like the single index case.
However, instead of 
(\ref{eq:2indRFOU_0050}),
a different RFOU with two indices that has a simpler spectral density will be considered here.

Let $\mathbf{D}_t^\alpha = \bigl(-d^2/dt^2\bigr)^{\alpha/2} $  ,  $\alpha > 0$ 
be the one-dimensional Riesz derivative defined by 
\begin{align}
  \mathbf{D}_t^\alpha{f} & = \bigl(-d^2/dt^2\bigr)^{\alpha/2}{f} = \mathcal{F}^{-1}\bigl(|k|^\alpha \widetilde{f}(k)\bigr), 
\label{eq:2indRFOU_0060}
\end{align}
where $\widetilde{f}(k)$ is the Fourier transform of $f(t)$.
By replacing the Weyl fractional derivative in 
\mbox{(\ref{eq:2indRFOU_0010})}
with the Riesz derivative results in the fractional Langevin equation of Reisz type for FOU with two indices $Y_{\alpha\beta,\lambda}(y)$:

\begin{align}
    \bigl(\mathbf{D}_t^{2\beta} + \lambda^{2\beta}\bigr)^{\alpha/2} Y_{\alpha\beta,\lambda}(t) & = \eta(t).
\label{eq:2indRFOU_0070}
\end{align}
The solution of 
(\ref{eq:2indRFOU_0070})
is given by
\begin{align}
     Y_{\alpha\beta,\lambda}(t) & = \frac{1}{\sqrt{2\pi}}\int_{-\infty}^\infty \frac{e^{ikt}\widetilde{\eta}(k)}{\bigl(|k|^{2\beta} + \lambda^{2\beta}\bigr)^\alpha} dk.
\label{eq:2indRFOU_0080}
\end{align}
The spectral density of this process has a simpler form as compared with 
(\ref{eq:2indRFOU_0040}):
\begin{align}
   S_{\alpha\beta}(k) & = \frac{1}{2\pi\bigl(|k|^{2\beta} + \lambda^{2\beta}\bigr)^\alpha}.
\label{eq:2indRFOU_0090}
\end{align}
The covariance function $C_{\alpha\beta}(t)$  can be obtained by taking the Fourier transform of $S_{\alpha\beta}(k)$. 
However, it does not in general has a closed analytic form. 
The variance is
\begin{align}
  C_{\alpha\beta,\lambda}(0) & = \frac{1}{\pi} \int_0^\infty \frac{1}{\bigl(|k|^{2\beta} + \lambda^{2\beta}\bigr)^\alpha} dk
                          = \frac{\Gamma(1/2\beta)\Gamma(\alpha-1/2\beta)}{2\pi\beta\Gamma(\alpha)}\lambda^{1-2\alpha\beta},     
\label{eq:2indRFOU_0100}
\end{align}
again the condition $\alpha\beta > 1/2$   is imposed to ensure the variance is finite.

Note that both $X_{\alpha\beta,\lambda}(t)$  and $Y_{\alpha\beta,\lambda}(t)$  can be regarded as two different generalisations of FOU of single index to two indices. 
These two processes have similar long and short time asymptotic properties 
\cite{LimLiTeo2008,LimTeo2009c}.
However, only the reduced process associated with  $Y_{\alpha\beta,\lambda}(t)$  will be considered here.

The reduced process associated with $Y_{\alpha\beta,\lambda}(t)$, again denoted by $B_{\alpha\beta,\lambda}(t)$, 
has the following spectral representation:
\begin{align}
  B_{\alpha\beta,\lambda}(t) & = Y_{\alpha\beta,\lambda}(t) - Y_{\alpha\beta,\lambda}(0)
   = \frac{1}{\sqrt{2\pi}} \int_{-\infty}^\infty \frac{\bigl(e^{ikt} -1\bigr)\widetilde{\eta}(k)}{\bigl(|k|^{2\beta} + \lambda^{2\beta}\bigr)^{\alpha/2}} dk .
\label{eq:2indRFOU_0110}
\end{align}
Its covariance is given by
\begin{align}
  \breve{C}_{\alpha\beta,\lambda}(t,s) & = \Bigl<B_{\alpha\beta,\lambda}(t)B_{\alpha\beta,\lambda}(s)\Bigr>
                                          = \frac{1}{2\pi} \int_{-\infty}^\infty \frac{\bigl(e^{ik|t - s|} - e^{ik|t|} - e^{ik|s|} +1\bigr)}{\bigl(|k|^{2\beta} + \lambda^{2\beta}\bigr)^\alpha} dk,
\label{eq:2indRFOU_0130}
\end{align}
and its variance is
\begin{align}
  \breve{\sigma}_{\alpha\beta,\lambda}^2 (t) & = \breve{C}_{\alpha\beta,\lambda}(t,t) 
                                           = \frac{1}{\pi}\int_{\mathbb{R}} \frac{1 - e^{ikt}}{\bigl(|k|^{2\beta} + \lambda^{2\beta}\bigr)^\alpha}dk.
\label{eq:2indRFOU_0140}
\end{align}
The covariance 
can be expressed in a form similar to TFBM with single index.
\begin{align}
  \breve{C}_{\alpha\beta,\lambda}(t,s) & = \frac{1}{2}\Bigl[c_t|t|^{2H} + c_s|s|^{2H} - c_{t-s}|t-s|^{2H}  \Bigr],
\end{align}
where $2\alpha\beta-1 =2H$, and
\begin{align}
  c_t & = \frac{\Gamma(1/2\beta)\Gamma(H/\beta)}{\pi\beta\Gamma\bigl[(2H+1)/2\beta\bigr]\lambda^{2H}}
        - \frac{1}{\pi|t|^{2H}}\int_{-\infty}^\infty \frac{e^{ikt}}{\bigl(|k|^{2\beta} + \lambda^{2\beta}\bigr)^{(2H+1)/2\beta}} dk.
\end{align}

As it will be shown below that $B_{\alpha\beta,\lambda}(t)$  satisfies the same properties as TFBM $B_{\alpha,\lambda}(t)$, 
so it can be regarded as TFBM with two indices.

\subsection{Scaling property}
\noindent
First, we verify the scaling property for $Y_{\alpha\beta,\lambda}(t)$
\begin{align}
  \Bigl<Y_{\alpha\beta,\lambda}(r(t+\tau))Y_{\alpha\beta,\lambda}(rt)\Bigr> & = \frac{1}{\pi}\int_0^\infty \frac{\cos(k|r\tau|)}{\bigl(|k|^{2\beta} + \lambda^{2\beta}\bigr)^\alpha}dk \nonumber \\
     &  = \frac{|r\tau|^{2\alpha\beta-1}}{\pi}\int_0^\infty \frac{\cos(k)}{\bigl(|k|^{2\beta} + (\lambda|r\tau|)^{2\beta}\bigr)^\alpha}dk  \nonumber \\
     &  = \frac{r^{2\alpha\beta-1}|\tau|^{2\alpha\beta-1}}{\pi}\int_0^\infty \frac{\cos(k)}{\bigl(|k|^{2\beta} + (\lambda|r\tau|)^{2\beta}\bigr)^\alpha}dk \nonumber \\
     & =   r^{2\alpha\beta-1}\Bigl<Y_{\alpha\beta,r\lambda}((t+\tau))Y_{\alpha\beta,r\lambda}(t)\Bigr>,
\label{eq:2indRFOU_0150}
\end{align}
where $Y_{\alpha\beta,r\lambda}(t)$ is the same process as $Y_{\alpha\beta,\lambda}(t)$ with $\lambda$ replaced by $r\lambda$. 
By expressing $\breve{C}_{\alpha\beta,\lambda}(t+\tau,t)$ in terms of the covariance of $Y_{\alpha\beta,\lambda}(t)$, one has 
\begin{align}
\breve{C}_{\alpha\beta,\lambda}(t,t+\tau) & 
                                            = C_{\alpha\beta,\lambda}(\tau) - C_{\alpha\beta,\lambda}(t+\tau)
                                              - C_{\alpha\beta,\lambda}(t) + C_{\alpha\beta,\lambda}(0).
\label{eq:2indRFOU_0170}
\end{align}
With the help of
(\ref{eq:2indRFOU_0150}),
one obtains the scaling property for $B_{\alpha\beta,\lambda}$:
\begin{align}
  \Bigl<B_{\alpha\beta,\lambda}(rt)B_{\alpha\beta,\lambda}(rs)\Bigr> & = r^{2\alpha\beta-1}\Bigl<B_{\alpha\beta,r\lambda}(t)B_{\alpha\beta,r\lambda}(s)\Bigr>.
\label{eq:2indRFOU_0180}
\end{align}

\subsection{Stationary increments}
\noindent
The increment process of $B_{\alpha\beta,\lambda}$ is stationary which follows from the stationarity of $Y_{\alpha\beta,\lambda}$ by noting that
\begin{align}
\Delta_\tau B_{\alpha\beta,\lambda}(t) & =  B_{\alpha\beta,\lambda}(t+\tau) - B_{\alpha\beta,\lambda}(t) = Y_{\alpha\beta,\lambda}(t+\tau) - Y_{\alpha\beta,\lambda}(t).
\label{eq:2indRFOU_0190}
\end{align}
Its covariance is
\begin{align}
  \Bigl<\Delta_\tau B_{\alpha\beta,\lambda}(t_1)\Delta_\tau B_{\alpha\beta,\lambda}(t_2)\Bigr>
           = \frac{1}{2\pi} \int_{-\infty}^\infty \frac{2e^{ik|t_2-t_1|} - e^{ik|t_2-t_1+\tau|} - e^{ik|t_2-t_1-\tau|}}{\bigl(|k|^{2\beta} + \lambda^{2\beta}\bigr)^\alpha} dk.
\label{eq:2indRFOU_0200}
\end{align}

\subsection{Locally self-similar property}
\noindent
One can verify that $B_{\alpha\beta,\lambda}(t)$ is locally asymptotically self-similar of order $\alpha\beta -1/2$   by showing 
\begin{align}
   \lim_{\epsilon \to 0} \Biggl[\frac{B_{\alpha\beta,\lambda}(t_\circ+\epsilon{u}) - B_{\alpha\beta,\lambda}(t_\circ)}{\epsilon^{\alpha\beta-1}}\Biggr]
    & \ \hat{=} \ T_{t_\circ}(u),
\label{eq:2indRFOU_0270}
\end{align}
where the tangent process $T_{t_\circ}(u)$  at point $t_\circ$  is a FBM indexed by Hurst index $H = \alpha{\beta{ - 1/2}}$. 
Since the increment process
\begin{align}
  \Delta_{\epsilon{u}}B_{\alpha\beta,\lambda}(t_\circ) & = B_{\alpha\beta,\lambda}(t_\circ + \epsilon{u})  - B_{\alpha\beta,\lambda}(t_\circ) 
                         = Y_{\alpha\beta,\lambda}(t_\circ+\epsilon{u}) - Y_{\alpha\beta,\lambda}(t_\circ)
                         = \Delta_{\epsilon{u}} Y_{\alpha\beta,\lambda}(t_\circ)
\label{eq:2indRFOU_0280}
\end{align}
so that instead of 
(\ref{eq:2indRFOU_0270}),
one can use
\begin{align}
  \lim_{\epsilon \to 0} \Biggl[\frac{Y_{\alpha\beta,\lambda}(t_\circ+\epsilon{u}) - Y_{\alpha\beta,\lambda}(t_\circ)}{\epsilon^{\alpha\beta-1}}\Biggr]
    & \ \hat{=} \ T_{t_\circ}(u),
\label{eq:2indRFOU_0280}
\end{align}
for the verification of locally asymptotically self-similarity.

Before we can verify 
(\ref{eq:2indRFOU_0270})
or 
(\ref{eq:2indRFOU_0280})
we need to consider the leading term of the variance of the increment process
\begin{align}
  \breve{\sigma}_{\alpha\beta,\Delta}^2(t) & = \sigma_{\alpha\beta,\Delta}^2(t)
                                     = \Bigl<\bigl(Y_{\alpha\beta,\lambda}(t+s) - Y_{\alpha\beta,\lambda}(s)\bigr)^2\Bigr>
     = 2\bigl(C_{\alpha\beta}(0) - C_{\alpha\beta}(t)\bigr)
\label{eq:2indRFOU_0290}
\end{align}
For $1/2 < \alpha\beta < 3/2$,
\begin{subequations}
\label{eq:2indRFOU_0300}
\begin{align}
  \sigma_{\alpha\beta,\Delta}^2(t) & =\frac{4}{\pi}  \int_0^\infty \frac{\sin^2\bigl(k|t|/2\bigr)}{\bigl(k^{2\beta} + \lambda^{2\beta}\bigr)^\alpha}dk
                                     =\frac{4|t|^{2\alpha\beta-1}}{\pi} \int_0^\infty  \frac{\sin^2(k/2)}{\bigl(k^{2\beta} + \lambda^{2\beta}|t|^{2\beta}\bigr)^\alpha}dk \nonumber \\
                                   & =\frac{4|t|^{2\alpha\beta-1}}{\pi} \int_0^\infty k^{-\alpha\beta} \sin^2(k/2) dk + o\bigl(|t|^{2\alpha\beta-1}\bigr) \nonumber \\
                                   & = \frac{|t|^{2\alpha\beta-1}}{\Gamma(\alpha\beta)\cos(\alpha\beta\pi)}  + o\bigl(|t|^{2\alpha\beta-1}\bigr), &
                                   & \text{as} \ t \to 0.
\label{eq:2indRFOU_0300a}
\end{align}
Or with $2\alpha\beta - 1 = 2H$,
\begin{align}
  \sigma_{\alpha\beta,\Delta}^2(t)  & \sim -\frac{|t|^{2H}}{\Gamma(\alpha\beta)\cos(\alpha\beta\pi)}  + o\bigl(|t|^{2H}\bigr) , &
     & \text{as} \ t \to 0,
\end{align}
\end{subequations}
which shows that the short-time asymptotic behaviour of $\sigma_{\alpha\beta,\Delta}^2(t)$ 
varies as $t^{2\alpha\beta-1}$ or $t^{2H}$.
Note that \mbox{(\ref{eq:2indRFOU_0300})} can be re-expressed in the same form as 
(\ref{eq:fracOU_0120})
by using 
$\breve{\sigma}_{\alpha\beta,\Delta}(t) = 2C\bigl(C_{\alpha\beta}(0) - C_{\alpha\beta}(t)\bigr)$,
hence
$B_{\alpha\beta,\lambda}(t)$ is locally self-similar.

Now, consider the covariance of the tangent process at the point $t_\circ$:
\begin{gather} 
   \lim_{\epsilon \to 0} \Biggl[
                \Biggl<
                  \bigg(
                    \frac{B_{\alpha\beta,\lambda}(t_\circ+\epsilon{u}) - B_{\alpha\beta,\lambda}(t_\circ)}{\epsilon^{\alpha\beta-1}}
                  \biggr)
                  \biggl(
                    \frac{B_{\alpha\beta,\lambda}(t_\circ+\epsilon{v}) - B_{\alpha\beta,\lambda}(t_\circ)}{\epsilon^{\alpha\beta-1}}
                  \biggl)
                \Biggr>
                \Biggr] \qquad\qquad\qquad \nonumber \\
  \qquad\qquad\qquad\qquad     = \lim_{\epsilon \to 0} \Biggl[
                                              \frac{\sigma_{\alpha\beta,\Delta}^2(\epsilon{u})
                                               + \sigma_{\alpha\beta,\Delta}^2(\epsilon{v})
                                               - \sigma_{\alpha\beta,\Delta}^2(\epsilon(u-v))
                                               }{2\epsilon^{\alpha\beta-1}}
                                             \Biggr] \nonumber \\
  \qquad\qquad\qquad\qquad\qquad\qquad     =  \frac{1}{\Gamma(\alpha\beta)\cos(\alpha\beta\pi)}
                                               \bigl(|u|^{2\alpha\beta-1} + |v|^{2\alpha\beta-1} - |u-v|^{2\alpha\beta-1}\bigr) ,
\label{eq:2indRFOU_0310}
\end{gather}
which is just the covariance  
$\Bigl<B_H(u)B_H(v)\Bigr>$
of the fractional Brownian motion if we identify $H$ with $2\alpha\beta -1$.

Note that the spectral density of $Y_{\alpha\beta,\lambda}(t)$ given by 
(\ref{eq:2indRFOU_0090})
has the same functional form as the Linnik probability density function. 
The analytic properties of the latter depend on the arithmetic nature of both parameters $\alpha$ and $\beta$; 
and the conditions imposed on $\alpha$ and $\beta$ are rather complicated and are not of practical interest 
\cite{KOH1995ab,ErdoganOstrovskii1998,LimTeo2010,GradshteynRyzhik2000}.
Therefore, different methods have been adopted here for studying the asymptotic behaviour of the covariance function $C_{\alpha\beta}(t)$
that are sufficient for most practical purposes. 
In the above discussion, the values of  $\alpha$ and $\beta$
have been confined to  $1/2 < \alpha\beta < 3/2$, which correspond  to  $0 < H < 1$ with $\alpha\beta = H +1/2$.

Note that in the small time asymptotic behaviour 
(\ref{eq:2indRFOU_0300}),
$\sigma_{\alpha\beta,\Delta}^2(t)$ depends on $\alpha$ and $\beta$ through the product $\alpha\beta$.
Let $\alpha = \alpha^\prime/\beta$ such that $\alpha\beta = \alpha^\prime$.
Then $\breve{\sigma}_{\alpha\beta,\Delta}^2(t)$ depends on $\alpha^\prime$ and independent of $\beta$ as $t \to 0$.
On the other hand, the long-time asymptotic behaviour of the covariance varies as $t^{-(1+2\beta)}$
which is independent of $\alpha$.
Thus, in contrast to FBM and TFBM, 
the short-time property such as fractal dimension, and the long-time behaviour like long-range dependence of TFBM with two indices 
$B_{\alpha\beta,\lambda}(t)$
can be separately characterized by using two different parameters.

\subsection{Long range dependence}
\noindent
\label{sec:LongRange}
First, one note that $Y_{\alpha\beta,\lambda}(t)$ is a short memory process. 
Its covariance can be expressed as
\begin{align}
  C_{\alpha\beta,\lambda}(\tau) & = \frac{1}{\pi} \int_0^\infty \frac{\cos(k|\tau|)}{\bigl(|k|^{2\beta} + \lambda^{2\beta}\bigr)^\alpha} dk
                              = \frac{1}{\pi} \text{Im} \int_0^\infty \frac{e^{-ut}}{\bigl(e^{-i\beta\pi}u^{2\beta} + \lambda^{2\beta}\bigr)^\alpha}
\label{eq:2indRFOU_0210}
\end{align}
Using the series expansion
\begin{align*}
  \frac{1}{(1+z)^\alpha} & = \sum_{j=0}^\infty \frac{(-1)^j\Gamma(\alpha+j)}{j!\Gamma(\alpha)} z^j, 
\end{align*}
and substitution $u$ by $u/\tau$, one gets for $0 < \beta < 1$ the following $\tau \to \infty$ asymptotic expression
\begin{align}
  C_{\alpha\beta,\lambda}(\tau) & = \frac{1}{\pi}\text{Im}\Biggl[\frac{1}{\tau}\int_0^\infty \frac{e^{-u} du}{\bigl(e^{-i\beta\pi}(u/\tau)^{2\beta} + \lambda^{2\beta}\bigr)^\alpha}\Biggr] \nonumber \\
         & = \frac{1}{\pi}\text{Im}\Biggl[\frac{1}{\tau}\int_0^\infty e^{-u}
               \sum_{j=0}^\infty  \frac{(-1)^j\Gamma(\alpha+j)}{j!\Gamma(\alpha)} 
                 e^{-i\beta{j}\pi} \lambda^{-2\beta(\alpha+j)} (u/\tau)^{2\beta{j}} du
           \Biggr]  \nonumber \\
         & = \frac{1}{\pi\Gamma(\alpha)}
              \sum_{j=1}^\infty \frac{(-1)^{j+1} \lambda^{-2\beta(\alpha+j)}\Gamma(\alpha+j)\Gamma(1+2\beta{j})\sin(\beta{j}\pi)}{j!}
               \tau^{-(2\beta{j}+1)} .
\label{eq:2indRFOU_0220}
\end{align}
Thus, the leading term in the limit is 
$\frac{\alpha}{\pi}\lambda^{-2\beta(1+\alpha)}\Gamma(1+2\beta)\sin(\beta\pi)\tau^{-(1+2\beta)}$,
which gives polynomial decay $\tau^{-(1+2\beta)}$ for the covariance $C_{\alpha\beta}(\tau)$ as $\tau \to \infty$.
This implies $Y_{\alpha\beta,\lambda}(t)$ is a short memory process.

On the other hand, $B_{\alpha\beta,\lambda}(t)$ is long memory process. 
One can use the same argument as for $B_{\alpha,\lambda}(t)$  to show that $B_{\alpha\beta,\lambda}(t)$  is LRD.
Note that the $\tau \to \infty$ limit of the correlation function of $B_{\alpha\beta,\lambda}(t)$ satisfies
\begin{align}
  \breve{R}_{\alpha\beta,\lambda}(t,t+\tau) & \sim  \frac{1}{2}\sqrt{\frac{\breve{\sigma}_{\alpha\beta,\lambda}^2(t)}{\sigma_{\alpha\beta,\lambda}^2(t)}}  > 0,
\label{eq:2indRFOU_0230}
\end{align}
where $\sigma_{\alpha\beta,\lambda}(t)$ is just a constant given by  
(\ref{eq:2indRFOU_0100}).
This implies $B_{\alpha\beta,\lambda}(t)$ is a LRD process.

\subsection{Fractal dimension}
\noindent
Fractal dimension is a local property and since the process $B_{\alpha\beta,\lambda}(t)$ satisfies locally asymptotically self-similar property, 
it behaves locally like fractional Brownian motion. 
One expects the fractal dimension of the graph of $B_{\alpha\beta,\lambda}(t)$ is the same as that of FBM.

$B_{\alpha\beta,\lambda}(t)$ satisfies
\begin{align}
  \wideparen{\sigma}_{\alpha\beta,\Delta}(t) & = \Bigl<\bigl(B_{\alpha\beta,\lambda}(t+s) - B_{\alpha\beta,\lambda}(s)\bigr)^2\Bigr> \nonumber \\
                                         & = \Bigl<\bigl(Y_{\alpha\beta,\lambda}(t+s) - Y_{\alpha\beta,\lambda}(s)\bigr)^2\Bigr>      
                                           = \sigma^2_{\alpha\beta,\Delta}(t) \leq A|\tau|^{2\alpha\beta}.
\label{eq:2indRFOU_0320}
\end{align}
Thus, almost surely the sample path of $B_{\alpha\beta,\lambda}(t)$ is H\"{o}lderian of order $(\alpha\beta -1/2) - \epsilon$,
for all $\epsilon > 0$.
Hence the fractal dimension of the graph of $B_{\alpha\beta,\lambda}(t)$ is a.e. equal to $5/2 - \alpha\beta$ 
(or $2 - H$) for $1/2 < \alpha\beta < 3/2 \ (0 < H <1)$

Note that in the small time asymptotic behaviour given by 
\mbox{(\ref{eq:2indRFOU_0300})}
and 
\mbox{(\ref{eq:2indRFOU_0320})}
depends on 
$t^{2\alpha\beta-1}$.
Let $\alpha = \alpha^\prime/\beta$ such that $\alpha\beta =\alpha^\prime$.
One then has that local properties such as locally asymptotically self-similarity and fractal dimension depend on $\alpha^\prime$   
and independent of $\beta$. 
On the other hand, the long-time asymptotic behaviour of the covariance varies as $t^{-(1+2\beta)}$,
which is independent of $\alpha$.
In contrast to FBM and TFBM, it is possible to separately characterize the short-time property such as fractal dimension, 
and the long-time behaviour like long-range dependence of $B_{\alpha\beta,\lambda}(t)$
by using two different parameters. 
The ability for a stochastic process to have separate characterization of fractal dimension and 
long-range dependency is a desirable property in the modeling of physical and geological phenomena.

Though it is not possible to express the covariance of $B_{\alpha\beta,\lambda}(t)$ n closed analytic form, 
one can show that it has the similar properties as TFBM  $B_{\alpha,\lambda}(t)$. 
Therefore, $B_{\alpha\beta,\lambda}(t)$ can be regarded as a generalization of TFBM from single index to two indices with a richer structure.


\section{Tempered Multifractional Brownian Motion}
\label{sec:TMfBM}
\noindent
TFBM with its stationary increments locally self-similarity provides a simple model for describing phenomena such as wind speed in 
Davenport's model. 
The main attractiveness of the model is its simplicity with each of the properties described by a single index 
$\alpha$ (or $H = \alpha -1/2$).
However, situation in the real world is more complicated, 
scaling property may vary with time (or position) and the system may have variable time-dependent memory. 
Therefore for a more realistic model it is necessary to introduce time-dependent or position-dependent scaling exponent. 
This can be achieved by 
extending TFBM to tempered multifractional Brownian motion (TMBM) with a variable index 
$\alpha(t)$,
in analogy with the generalization of FBM to MBM 
\mbox{\cite{Levy-VehelPeltier1996,BenassiJaffardRoux1997}}.
One can regard TMBM as the reduced multifractional Ornstein-Uhlenbeck process (RMOU) just like TFBM. 
Following the moving average definition of MBM 
\cite{Levy-VehelPeltier1996,BenassiJaffardRoux1997},
FOU of Weyl type given by 
(\ref{eq:fracOU_0070})
can be generalised to multifractional Ornstein-Uhlenbeck process (MOU) by replacing the constant index $\alpha$  by a deterministic function $\alpha(t)$ 
\cite{LimTeo2007}:
\begin{align}
  X_{\alpha(t),\lambda}(t) & = \frac{1}{\Gamma\bigl(\alpha(t)\bigr)} \int_{-\infty}^t e^{-\lambda(t-u)} (t - u)^{\alpha(t)-1} \eta(u)du,
\label{eq:TMfBM_0010}
\end{align}
where variable index $\alpha(t)$ satisfies $\alpha(t) > 1/2$, is assumed to be H\"{o}lder continuous with
$\bigl|\alpha(t) - \alpha(s)\bigr| \leq k|t - s|^\beta$, $k > 0$, $\beta > 0$.
For $s < t$, the covariance of $X_{\alpha(t),\lambda}(t)$ is given by
\begin{align}
  \Bigl<X_{\alpha(t),\lambda}(t)X_{\alpha(s),\lambda}(s)\Bigr> & = \frac{e^{-\lambda(t+s)}}{\Gamma\bigl(\alpha(t)\bigr)\Gamma\bigl(\alpha(s)\bigr)} \int_{-\infty}^{\min(t,s)} (t -u)^{\alpha(t) -1} (s -u)^{\alpha(s) -1} e^{2\lambda{u}} du \nonumber \\
& = \frac{e^{-\lambda(t-s)}}{\Gamma\bigl(\alpha(t)\bigr)\Gamma\bigl(\alpha(s)\bigr)} \int_{-\infty}^s (u)^{\alpha(s) -1} (u+t- s)^{\alpha(t) -1} e^{-2\lambda{u}} du \nonumber\\
  & =  \frac{e^{-\lambda(t-s)}(t - s)^{2\alpha_{+}(t,s)-1}}{\Gamma\bigl(\alpha(t)\bigr)}
       \Psi\bigl(\alpha(s),2\alpha_{+}(s,t), 2\lambda(t-s)\bigr) ,
\label{eq:TMfBM_0020}
\end{align}
where  $\Psi(\alpha,\gamma;z)$ is the confluent hypergeometric function, which is also known as Kummer 
function and is denoted by $U(\alpha,\gamma,z)$, and 
$\alpha_{+}(s,t)  = \bigl(\alpha(t) + \alpha(s)\bigr)/2$ 
(3.383 of \cite{GradshteynRyzhik2000}
has been used).
Another possible (equivalent) way of defining MFOU is based on the spectral representation 
(\ref{eq:fracOU_0110})
of FOU:
\begin{align}
  X_{\alpha(t),\lambda}(t) & = \frac{1}{\sqrt{2\pi}}\int_{-\infty}^\infty \frac{e^{-ikt}\widetilde{\eta}(k)dk}{\bigl(-ik + \lambda\bigr)^{\alpha(t)}} .
\label{eq:TMfBM_0030}
\end{align}
The covariance is 
\begin{align}
  \Bigl<X_{\alpha(t),\lambda}(t)X_{\alpha(s),\lambda}(s)\Bigr> & = \frac{1}{2\pi}\int_{-\infty}^\infty \frac{e^{-ik(t-s)}dk}{\bigl(-ik + \lambda\bigr)^{\alpha(t)}\bigl(ik + \lambda\bigr)^{\alpha(s)}}  \nonumber\\
    & = \frac{(t-s)^{\alpha_{+}(s,t)-1}}{\Gamma\bigl(\alpha(t)\bigr)(2\lambda)^{\alpha_{+}(s,t)}}
        W_{\alpha_{-}(s,t), 1/2 - \alpha_{+}(s,t) }\bigl(2\lambda(t-s)\bigr) ,
\label{eq:TMfBM_0040}
\end{align}
where $W(\cdot)$ is the Whittaker function, and $\alpha_{-}(t,s) = \bigl(\alpha(t) - \alpha(s)\bigr)/2$
 (see 3.384 of 
\cite{GradshteynRyzhik2000}
).

There are two definitions of the MOU
$X_{\alpha(t),\lambda}(t)$
given by 
\mbox{(\ref{eq:TMfBM_0010}) and (\ref{eq:TMfBM_0030})}.
For the proof of the equivalence of these two representations of a centred Gaussian process suffice to show that the covariance functions given by 
\mbox{(\ref{eq:TMfBM_0020}) and (\ref{eq:TMfBM_0040})}
are the same 
(see \ref{sec:FWequivalent}
for the verification). 
Note that just like in the fractional case, the variance and covariance functions 
of TMBM are divergent as $\lambda \to 0$ \mbox{\cite{LimTeo2007}}.

Instead of defining TMBM based on the reduced process of $X_{\alpha(t),\lambda}(t)$,
another version of MOU will be used for this purpose. 
The reason for introducing a different MOU is that it leads to a TMBM with 
a covariance function which has the same functional form as the covariance of TFBM with the constant index replaced by variable index. 
By letting $\beta=1$ and replacing $\alpha$ by $\alpha(t)$ in 
(\ref{eq:2indRFOU_0080})
leads to
\begin{align}
  Y_{\alpha(t),\lambda}(t) & = \frac{1}{\sqrt{2\pi}}\int_{-\infty}^\infty \frac{e^{ikt}\widetilde{\eta}(k)dk}{\bigl(|k|^2 + \lambda^2\bigr)^{\alpha(t)/2}} .
\label{eq:TMfBM_0060}
\end{align}
Here we have use $Y_{\alpha(t),\lambda}(t)$ instead of $X_{\alpha(t),\lambda}(t)$ to distinguish from the previous version of MOU. 
Ruiz-Medina \textit{et al}.\@ \cite{Ruiz-MedinaAnguloAnh2003} 
claimed that it is the mean square solution to the one-dimensional fractional Bessel equation with variable order
\begin{align}
  \bigl(\mathbf{D}_t^2 + \lambda^2\bigr)^{\alpha(t)/2} Y_{\alpha(t),\lambda}(t) & = \eta(t),
\label{eq:TMfBM_0050}
\end{align}
where $\mathbf{D}_t^{\alpha(t)} = \bigl(-d^2/dt^2\bigr)^{\alpha(t)/2}$ is the one-dimensional Riesz derivative of variable order
$\alpha(t)$, with $\alpha(t) > 1/2$. 
The covariance of $Y_{\alpha(t),\lambda}(t)$ is given by
\begin{align}
  \Bigl<Y_{\alpha(t),\lambda}(t)Y_{\alpha(s),\lambda}(s)\Bigr> & = \frac{1}{\sqrt{\pi}\Gamma\bigl(\alpha_{+}(s,t)\bigr)}
                                                           \left(\frac{|t-s|}{2\lambda}\right)^{\alpha_{+}(s,t) - 1/2}
                                                           K_{\alpha_{+}(s,t) - 1/2}\bigl(\lambda|t - s|\bigr) ,     
\label{eq:TMfBM_0070}
\end{align}
and the variance is 
\begin{align}
  \Bigl<\bigl(Y_{\alpha(t),\lambda}(t)\bigr)^2\Bigr>  & = \frac{\Gamma\bigl(2\alpha(t)-1\bigr)}{\Bigl(\Gamma\bigl(\alpha(t)\bigr)\Bigr)^2 \bigl(2\lambda\bigr)^{2\alpha(t)-1}} .
\label{eq:TMfBM_0071}
\end{align}
Note that 
(\ref{eq:TMfBM_0050})
can also be regarded as a special case of the fractional Riesz-Bessel process of variable order 
\cite{LimTeo2008a}.
The process defined by 
(\ref{eq:TMfBM_0060})
with covariance 
(\ref{eq:TMfBM_0070})
appears to be the appropriate MOU since its reduced process leads to TMBM as generalization of TFBM in line with the extension of FBM to MBM. 
Consider the RMOU
\begin{align}
  B_{\alpha(t),\lambda}(t) & = Y_{\alpha(t),\lambda}(t) - Y_{\alpha(t),\lambda}(0),
\label{eq:TMfBM_0080}
\end{align}
with covariance given by
\begin{align}
  \breve{C}_{\alpha(t),\lambda}(t - s) & =  \Bigl<B_{\alpha(t),\lambda}(t)B_{\alpha(s),\lambda}(s)\Bigr> \nonumber \\
                                       & = \frac{1}{\sqrt{\pi}\alpha_{+}(s,t)}
                                            \Biggl[
                                            \left(\frac{|t-s|}{2\lambda}\right)^{\alpha_{+}(s,t) - 1/2}
                                            K_{\alpha_{+}(s,t) - 1/2}\bigl(\lambda|t - s|\bigr) \nonumber \\
                                      & \qquad\qquad   - \left(\frac{|t|}{2\lambda}\right)^{\alpha_{+}(s,t) - 1/2}
                                           K_{\alpha_{+}(s,t) - 1/2}\bigl(\lambda|t|\bigr)   \nonumber \\
                                      & \qquad\qquad\qquad 
                                        - \left(\frac{|s|}{2\lambda}\right)^{\alpha_{+}(s,t) - 1/2}
                                           K_{\alpha_{+}(s,t) - 1/2}\bigl(\lambda|s|\bigr)  
                                           \Biggr]  \nonumber \\
                                      & \qquad\qquad\qquad\qquad 
                                        + \frac{\Gamma(2\alpha_{+}(s,t) - 1)}{\bigl(\Gamma(\alpha_{+}(s,t))\bigr)^2(2\lambda)^{2\alpha_{+}(s,t) - 1}}.
\end{align}
By letting $\alpha_{+}(s,t) - 1/2 = H_{+}(t,s) = \bigl(H(s) + H(t)\bigr)/2$, 
the covariance of TMFM can be expressed in the following form
\begin{subequations}
\label{eq:TMfBM_0090}
\begin{align}
  \breve{C}_{\alpha_{+}(s,t),\lambda}(t - s) & = \frac{1}{2} \Bigl[
                                          c_t\bigl(H_{+}(s,t)\bigr)|t|^{2H_{+}(s,t)}
                                         + c_s\bigl(H_{+}(s,t)\bigr)|s|^{2H_{+}(s,t)} \nonumber \\
                                    &  \qquad\qquad\qquad\qquad\qquad
                                         - c_{t-s}\bigl(H_{+}(s,t)\bigr)|t-s|^{2H_{+}(s,t)}
                                         \Bigr],
\label{eq:TMfBM_0090a}
\end{align}
with
\begin{align}
  c_t\bigl(H_{+}(s,t)\bigr) & = \frac{2\Gamma\bigl(2H_{+}(s,t)\bigr)}{\Gamma\bigl(2H_{+}(s,t) + 1/2\bigr)(2\lambda)^{2H_{+}(s,t)}}
                             \nonumber \\
                           & \qquad\qquad\qquad
                             - \frac{2}{\sqrt{\pi}\Gamma\bigl(H_{+}(s,t)+1/2\bigr)}
                               \left(\frac{1}{2\lambda|t|}\right)^{H_{+}(s,t)} K_{H_{+}(s,t)}\bigl(\lambda|t|\bigr).
\label{eq:TMfBM_0090b}
\end{align}
\end{subequations}
The variance is given by
\begin{align}
  \breve{\sigma}_{\alpha(t),\lambda}(t) & =  \Bigl<\bigl(B_{\alpha(t),\lambda}\bigr)^2\Bigl> \nonumber \\
                                                        & = \frac{2}{\sqrt{\pi}\Gamma\bigl(\alpha(t)\bigr)} 
                                                          \left[
                                                          \frac{\Gamma\bigl(2\alpha(t) -1\bigr)}
                                                               {\Bigl(\Gamma\bigl(\alpha(t)\bigr)\Bigr)^2(2\lambda)^{2\alpha(t)-1}}
                                                          - \left(\frac{|t|}{2\lambda}\right)^{\alpha(t) - 1/2}
                                                          K_{\alpha(t) - 1/2}\bigl(\lambda|t|\bigr)
                                                          \right].
\label{eq:TMfBM_0100}
\end{align}

The basic properties of TMBM are given below.

\subsection{Scaling property}
\noindent
Note that by naively replacing the index $\alpha$ by $\alpha(t)$ the scaling property (\ref{eq:fracOU_0150}) becomes
\begin{align}
  Y_{\alpha(t),\lambda}(bt) & \ \hat{=} \ b^{\alpha(t)-1/2} Y_{\alpha(t),b\lambda}(t) . 
\label{eq:TMfBM_0110}
\end{align}
The global scaling property 
(\ref{eq:TMfBM_0110})
fails to hold for TMBM, since one has
\begin{align}
  Y_{\alpha(t),\lambda}(s) & \ \hat{=} \ \left(\frac{s}{t}\right)^{\alpha(t)-1/2} Y_{\alpha(t),\lambda{s/t}}(t) . 
\label{eq:TMfBM_0111}
\end{align}
which implies the value of TMBM at time $s$ far apart from $t$ depends on the function $\alpha(t)$.
However, $Y_{\alpha(t),\lambda}(t)$ is locally asymptotically self-similar.

\subsection{Locally asymptotically self-similarity}
\noindent
TMBM is locally asymptotically self-similar of order $\alpha(t) - 1/2$ at a point $t$. 
In other words, the tangent process of $B_{\alpha(t),\lambda}(t)$ at a point $t$  is the FBM indexed by $\alpha(t) - 1/2$ (or $H(t)$).
For simplicity, assume  $\alpha(t)$ is H\"{o}lder continuous with $\bigl|\alpha(t) - \alpha(s)\bigr| \leq K|t - s|^\kappa$, and
$1/2 < \alpha(t) < \kappa + 1/2$ for all $t$.
Note that
\begin{align}
  B_{\alpha(t),\lambda}(t) - B_{\alpha(s),\lambda}(s) & = \bigl(Y_{\alpha(t),\lambda}(t) - Y_{\alpha(s),\lambda}(s)\bigr)  
                                                - \bigl(Y_{\alpha(t),\lambda}(0) - Y_{\alpha(s),\lambda}(0)\bigr) .         
\label{eq:TMfBM_0121}
\end{align}
Following a similar argument as for multifractional Riesz-Bessel process 
\cite{LimTeo2008a},
one has for $|t-s| \to 0$,
\begin{align}
  \bigl(Y_{\alpha(t),\lambda}(0) - Y_{\alpha(s),\lambda}(0)\bigr) & = O\bigl(|t-s|^{2\kappa}\bigr) ,          
\label{eq:TMfBM_0122}
\end{align}
and
\begin{align}
    \Bigl<\bigl(Y_{\alpha(t),\lambda}(t) - Y_{\alpha(s),\lambda}(s)\bigr)^2\Bigr> & = \frac{\Gamma\Bigl(\frac{1}{2} - \alpha(t)\Bigr)}
                                                                                    {2^{2\alpha(t) -1}\sqrt{\pi}\Gamma\bigl(\alpha(t)\bigr)}
                                                                                    |s - t|^{2\alpha(t)-1}
                                                                             + O\bigl(|s - t|^{2\alpha(t) - 1/2}\bigr) .     
\label{eq:TMfBM_0130}
\end{align}

Let
\begin{align}
  C_{\alpha(t)}(\epsilon;u,v) & = \Biggl<
                                \left(\frac{Y_{\alpha(t),\lambda}(t+\epsilon{u}) - Y_{\alpha(t),\lambda}(t)}{\epsilon^{\alpha(t) - 1/2}}\right)
                                \left(\frac{Y_{\alpha(t),\lambda}(t+\epsilon{v}) - Y_{\alpha(t),\lambda}(t)}{\epsilon^{\alpha(t) - 1/2}}\right)
                                \Biggr> \nonumber \\
     & = \frac{1}{2\epsilon^{2\alpha(t)-1}}
          \biggl[
              \Bigl<\bigl(Y_{\alpha(t),\lambda}(t+\epsilon{u}) - Y_{\alpha(t),\lambda}(t)\bigr)^2\Bigr> \nonumber \\
       & \qquad\qquad\qquad
         + \Bigl<\bigl(Y_{\alpha(t),\lambda}(t+\epsilon{v}) - Y_{\alpha(t),\lambda}(t)\bigr)^2\Bigr>  \nonumber \\
       & \qquad\qquad\qquad\qquad
         - \Bigl<\bigl(Y_{\alpha(t),\lambda}(t+\epsilon{v}) - Y_{\alpha(t),\lambda}(t+\epsilon{v})\bigr)^2\Bigr>
          \biggr].
\label{eq:TMfBM_0140}
\end{align}
For $\epsilon \to 0$, one has from 
(\ref{eq:TMfBM_0130})
\begin{align}
  \Bigl<\bigl(B_{\alpha(t),\lambda}(t+\epsilon{u}) - B_{\alpha(t),\lambda}(t))\bigr)^2\bigr> 
         & = \frac{\Gamma\bigl(1/2 - \alpha(t)\bigr)}{2^{2\alpha(t)-1}\sqrt{\pi}\Gamma\bigl(\alpha(t)\bigr)}
           \bigl(\epsilon|u|\bigl)^{2\alpha(t)-1} + O\bigl(\epsilon^{\gamma+ \alpha(t)-1/2}\bigr).
\label{eq:TMfBM_0150}
\end{align}
Therefore,
\begin{align}
  \lim_{\epsilon\to 0} C_{\alpha(t)}(\epsilon;u,v) & = \lim_{\epsilon\to 0} \frac{1}{2\epsilon^{2\alpha(t)-1}}
                                                   \Biggl[
                                                      \frac{\Gamma\bigl(1/2 - \alpha(t)\bigr)}{2^{2\alpha(t)-1}\sqrt{\pi}\Gamma\bigl(\alpha(t)\bigr)}
                                                     \Bigl(\bigl(\epsilon|u|\bigr)^{2\alpha(t) -1} + \bigl(\epsilon|v|\bigr)^{2\alpha(t) -1}\Bigr)
                                                 \nonumber \\
           & - \frac{\Gamma\bigl(1/2 - \alpha(t+\epsilon{v})\bigr)}{2^{2\alpha(t)-1}\sqrt{\pi}\Gamma\bigl(\alpha(t+\epsilon{v})\bigr)}
             \bigl(\epsilon|u-v|\bigr)^{2\alpha(t+\epsilon{v}) -1} + O\bigl(\epsilon^{\gamma+\min\bigl(\alpha(t),\alpha(t+\epsilon{v})\bigr)-1/2}\bigr)
             \Biggr]  \nonumber \\
     & = \frac{\Gamma\bigl(1/2 - \alpha(t)\bigr)}{2^{2\alpha(t)-1}\sqrt{\pi}\Gamma\bigl(\alpha(t)\bigr)}
       \Bigl(|u|^{2\alpha(t) -1} + |v|^{2\alpha(t) -1} - |u-v|^{2\alpha(t) -1}\Bigr), 
\label{eq:TMfBM_0160}
\end{align}
where H\"{o}lder continuity of $\alpha(t)$  implies
\begin{align}
  \frac{\Gamma\bigl(1/2 - \alpha(t+\epsilon{v})\bigr)}{2^{2\alpha(t)-1}\sqrt{\pi}\Gamma\bigl(\alpha(t+\epsilon{v})\bigr)}
  \bigl(\epsilon|u-v|\bigr)^{2\alpha(t+\epsilon{v}) -1}  
  &  \rightarrow  \frac{\Gamma\bigl(1/2 - \alpha(t)\bigr)}{2^{2\alpha(t)-1}\sqrt{\pi}\Gamma\bigl(\alpha(t)\bigr)} 
    \bigl(\epsilon|u-v|\bigr)^{2\alpha(t) -1}  .
\label{eq:TMfBM_0170}
\end{align}
(\ref{eq:TMfBM_0160})
is the covariance of MBM indexed by $\alpha(t) - 1/2$. 
Therefore, TMBM is locally asymptotically self-similar, and the tangent process at a point $t$ is the MBM with index 
$H(t) = \alpha(t) - 1/2$.

\subsection{Fractal dimension}
\noindent
With probability one, the Hausdorff dimension of the graph of the TMBM $B_{\alpha(t)}$ indexed by $\alpha(t)$ over the interval 
$I \in \mathbb{R}$ is $\frac{5}{2} - m_I\bigl[\alpha(t)\bigr]$, where $m_I\bigl[\alpha(t)\bigr] = \min\bigl\{\alpha(t); t \in I\bigr\}$.
First, note that the leading properties of the variances of increments of multifractional Brownian motion and TMBM are the same. 
Since only this leading behaviour is used to arrive at the local properties of the TMBM, 
one can easily infer that 
the local properties of the MBM holds verbatim for the TMBM
if we identify $H(t)$ with $\alpha(t) - 1/2$,
in particular the fractal dimension of the graph of $B_{\alpha(t),\lambda}(t)$
over an interval $I$ is $5/2 - \min[\alpha(t): t \in I]$
\cite{BenassiCohenIstas2003}.

\subsection{ Long range dependence}
\noindent
\noindent
One first shows that the non-stationary MOU $Y_{\alpha(t)}(t)$ is a short memory process. 
From the covariance 
(\ref{eq:TMfBM_0070})
and variance 
(\ref{eq:TMfBM_0071})
of  $Y_{\alpha(t)}(t)$, its correlation can be expressed as in the following form
\begin{align}
  R(t+\tau,t) & = g\bigl(\alpha(t+\tau),\alpha(t)\bigr)\left(\frac{|\tau|}{2\lambda}\right)^{\alpha_{+}(t,t+\tau)-1/2} K_{\alpha_{+}(t,t+\tau)-1/2}\bigl(\lambda|\tau|\bigr),
\end{align}
where $g\bigl(\alpha(t+\tau),\alpha(t)\bigr)$ is a funciton of $\alpha(t+\tau)$ and $\alpha(t)$, 
which is bounded above since $\alpha(\cdot)$ is bounded above. 
The exponential decay of $ K_{\alpha_{+}(t,t+\tau)-1/2}\bigl(\lambda|\tau|\bigr)$ thus ensures that MOU is SRD.

However, TMBM  $Y_{\alpha(t)}(t)$ is long-range dependent. 
Again, one can use a similar argument for TFBM. 
Note that whether the index of TFBM is a constant or time-dependent, 
there exists $t_\circ > 0$ such that for all $t > t_\circ$,
the variance of the process is positive; and moreover, 
it is bounded $\sigma^2(0) > \breve{\sigma}^2(t) > 0$.
By noting 
$\breve{\sigma}^2(t+\tau) \to \frac{2\Gamma\bigl(2\alpha(t)-1\bigr)}{\Gamma\bigl(\alpha(t)\bigr)^2(2\lambda)^{2\alpha(t)-1}} = 2\sigma^2(t)$
which is bounded above, one gets as $\tau \to \infty$

\begin{align}
  R(t,t+\tau) & \sim \frac{\breve{\sigma}_{\alpha(t),\lambda}^2(t)}{2\sqrt{\breve{\sigma}_{\alpha(t)}^2(t)\breve{\sigma}_{\alpha(t+\tau)}^2(t+\tau)}} 
                = \frac{1}{2}\sqrt{\frac{\breve{\sigma}_{\alpha(t),\lambda}^2(t)}{\sigma_{\alpha(t),\lambda}^2(t)}} > 0. 
\label{eq:TMfBM_0201}
\end{align}
Thus, TMBM is LRD.


\section{Concluding Remarks}
\label{sec:conclusion}
\noindent
TFBM and its properties are considered from the standpoint of RFOU.
The nice feature of this approach is that it facilitates and simplifies the study of the properties of TFBM. 
In contrast to FOU which is stationary and SRD, TFBM is a non-stationary LRD process. 
Despite of this, some of the properties of TFBM can be shown to be inherited from or direct consequences of the properties of FOU, 
which include the scaling behavior, 
spectral representation, 
and stationary increments of the process. 
In addition to the properties of TFBM studied previously by other authors 
\cite{MeerschaertSbzikar2013,MeerschaertSabzikar2014},
some different properties such as locally self-similarity and fractal dimension have also been considered.

Several generalizations of TFBM with some nice properties have been considered in this paper. 
For examples, TFBM with two indices allows separate characterization of fractal dimension and LRD property; 
and 
TMBM permits variable fractal dimension and memory. 
Note that TFBM and its various generalizations have the same local behavior as their corresponding FOU and MOU processes. 
In particular, TFBM and TMBM behave like FBM and MBM respectively in the small-time scales. 
However, their global properties are not the same, in particular, the tempered processes are long-range dependent, 
whereas the FOU with single index and two indices, 
and MOU are short memory processes. 
It is hoped that the various generalizations of TFBM considered in this paper can provide more flexibility and better modeling of physical, 
biological sciences, finance and other areas.



\appendix
\renewcommand{\theequation}{\Alph{section}\arabic{equation}}

\section{Tempered Fractional Gaussian Noise}
\label{sec:TFGN}
\noindent
Just like the case of fractional Gaussian noise (FGN) 
which can be regarded as the derivative of FBM in the sense of generalized functions, 
tempered fractional Gaussian noise (TFGN) can be defined in a similar way:
\begin{align}
  \xi_{\alpha,\lambda} & = \frac{dB_{\alpha,\lambda}(t)}{dt}.
\label{eq:TFGN_0010}
\end{align}
By using 
(\ref{eq:fracOU_0060})
one has
\begin{align}
  \xi_{\alpha,\lambda}(t) & =  D_t \Bigl[\bigl(e^{-\lambda{t}}I_t^{\alpha}e^{\lambda{t}}\bigr)\eta(t)\Bigr] \nonumber \\
               & = -\lambda e^{-\lambda{t}}I_t^{\alpha}e^{\lambda{t}}\eta(t) + e^{-\lambda{t}}D_tI_t^{\alpha}e^{\lambda{t}}\eta(t) \nonumber \\
               & = -\lambda X_{\alpha,\lambda}(t) + Z_{\alpha,\lambda},
\label{eq:TFGN_0020}
\end{align}
where
\begin{align}
  Z_{\alpha,\lambda}(t) & = e^{-\lambda{t}}D_tI_t^{\alpha}e^{\lambda{t}}\eta(t) = e^{-\lambda{t}}I_t^{\alpha-1}e^{\lambda{t}}\eta(t).
\label{eq:TFGN_0030}
\end{align}
Note that $\xi_{\alpha,\lambda}(t)$   is a combination of FOU and the fractional noise $Z_{\alpha,\lambda}(t)$.
\begin{align}
  C_\xi(t,s) & = C^{\alpha-1,\alpha-1}(t,s) - \lambda C^{\alpha-1,\alpha}(t,s) - \lambda C^{\alpha,\alpha-1}(t,s) +\lambda^2 C^{\alpha,\alpha}(t,s),
\label{eq:TFGN_0040}
\end{align}
with
\begin{align}
  C^{\mu,\nu}(t,s) & = \frac{1}{\Gamma(\mu)\Gamma(\nu)} \int_{-\infty}^t \int_{-\infty}^s
                       e^{-\lambda(t-u)}e^{-\lambda(s-v)} (t -u)^{\mu-1}(s -v)^{\nu-1} 
                     \bigl<\eta(u)\eta(v)\bigr>dudv  \nonumber \\
                 & = \frac{e^{-\lambda(t+s)}}{\Gamma(\mu)\Gamma(\nu)} \int_{-\infty}^{\min(t,s)}
                       e^{2\lambda{u}} (t -u)^{\mu-1}(s -u)^{\nu-1} du  
\label{eq:TFGN_0050}
\end{align}
Let $s < t$ then
\begin{align}
  C^{\mu,\nu}(t,s)  & = \frac{e^{-\lambda(t+s)}}{\Gamma(\mu)\Gamma(\nu)} \int_{-\infty}^s
                       e^{2\lambda{u}} (t -u)^{\mu-1}(s -u)^{\nu-1} du  \nonumber \\
                  & = \frac{e^{-\lambda(t-s)}}{\Gamma(\mu)\Gamma(\nu)} \int_0^\infty
                       e^{-2\lambda{u}} (t-s+u)^{\mu-1}u^{\nu-1} du \nonumber \\
                  & = \frac{e^{-\lambda(t-s)}(t-s)^{\mu+\nu-1}}{\Gamma(\mu)\Gamma(\nu)} \int_0^\infty
                       e^{-2\lambda(t-s)u} (1+u)^{\mu-1}u^{\nu-1} du \nonumber \\
                  & = \frac{1}{\Gamma(\mu)\bigl(2\lambda(t-s)^{(\mu+\nu)/2}\bigr)}W_{(\mu-\nu)/2,((\mu+\nu-1)/2}\bigl(2\lambda(t-s)\bigr),
\label{eq:TFGN_0060}
\end{align}
where $W_{\mu,\nu}(z)$ is the Whittaker function \cite[\#3.385.3]{GradshteynRyzhik2000}.

Note that the covariance (\ref{eq:TFGN_0060}) can also be expressed in terms of Kummer function $U(a,b,z)$ 
or confluent hypergeometric function $\Psi(a,b;z)$:
\begin{align}
  C^{\mu,\nu}(t,s)  & = \frac{e^{-\lambda(t-s)}(t-s)^{\mu+\nu-1}}{\Gamma(\mu)}U\bigl(\nu,\mu+\nu,2\lambda(t-s)\bigr).
\label{eq:TFGN_0070}
\end{align}
When $\mu=\nu$, one gets for $Re{z} > 0$
\begin{align}
  U(\nu+1/2,2\nu+1,z) = \frac{1}{\sqrt{\pi}}\frac{e^z}{(2z)^\nu}K_\nu(z).
\label{eq:TFGN_0080}
\end{align}
Thus
\begin{align}
  C^{\alpha,\alpha}(t,s) & = \frac{e^{-\lambda|t-s|}|t-s|^{2\alpha-1}}{\Gamma(\alpha)}U\bigl(\alpha,2\alpha;2\lambda|t-s|\bigr) \nonumber \\
                      & = \frac{1}{\sqrt{\pi}\Gamma(\alpha)}\left(\frac{|t-s|}{2\lambda}\right)^{a-1/2} K_{\alpha-1/2}\bigl(\lambda|t-s|\bigr),
\label{eq:TFGN_0090}
\end{align}
which is the covariance of FOU.

From the covariance 
(\ref{eq:TFGN_0040})
and 
(\ref{eq:TFGN_0060})
or 
(\ref{eq:TFGN_0070}),
one notes that the TFGN is a stationary process. 
The large time-lag $t-s = \tau \gg 0$ expansion
\begin{align}
  C^{\mu,\nu}(\tau) & = \frac{e^{\lambda\tau}\tau^{\mu-1}}{\Gamma(\mu)\Gamma(\nu)}
                      \int_0^\infty e^{-2\lambda{u}}\left(1+\frac{u}{\tau}\right)^{\mu-1}u^{\nu-1} du \nonumber \\
                  & = \frac{e^{\lambda\tau}\tau^{\mu-1}}{\Gamma(\mu)\Gamma(\nu)}
                       \sum_{n=0}^\infty \binom{\mu-1}{n} \int_0^\infty e^{-2\lambda{u}}\left(\frac{u}{\tau}\right)^nu^{\nu-1} du \nonumber \\
                  & = \frac{e^{\lambda\tau}\tau^{\mu-1}}{\Gamma(\mu)\Gamma(\nu)}
                       \sum_{n=0}^\infty \binom{\mu-1}{n} \frac{\Gamma(\nu+n)}{(2\lambda)^{\nu+n}\tau^n}.
\label{eq:TFGN_0100}
\end{align}
which shows that TFGN is SRD.

For small time lag, 
\begin{align}
   C^{\mu,\nu}(\tau) & = \frac{e^{-\lambda\tau}}{\Gamma(\mu)\Gamma(\nu)}
                       \sum_{n=0}^\infty \binom{\mu-1}{n} \int_0^\infty e^{2\lambda{u}}\tau^nu^{\mu-1-n}u^{\nu-1} du \nonumber \\
                   & = \frac{e^{-\lambda\tau}}{\Gamma(\mu)\Gamma(\nu)}
                       \sum_{n=0}^\infty \binom{\mu-1}{n} \tau^n \frac{\Gamma(\mu+\nu-n-1)}{(2\lambda)^{\mu+\nu-n-1}}
\end{align}
one has  $C^{\mu,\nu}(\tau) \sim \frac{\Gamma(\mu+\nu-1)}{\Gamma(\mu)\Gamma(\nu)}\frac{1}{(2\lambda)^{\mu+\nu-1}}$
as $\tau \to 0$.


\section{The Equivalent of the Two Representations of MOU}
\label{sec:FWequivalent}
To verify the equivalence of the two representations of MOU $X_{\alpha\beta,\lambda}(t)$ given by
(\ref{eq:TMfBM_0020}) and (\ref{eq:TMfBM_0040}),
one considers the relation between the confluent hypergeometric function and Whittaker function.

Consider the following function \# 13.4.4,  \cite{OlverLozierBoisvertEtAl2010}, page 326)
\begin{align}
  U(a,b,z) & = \frac{1}{\Gamma(a)}\int_0^\infty e^{-zt} t^{a-1} (1+t)^{b-a-1}dt,
           & \Re{a} > 0, |ph z| < \frac{1}{2}\pi.  
\label{eq:FWequivalent_0060}
\end{align}

Whittaker function in term of $U(a,b,z)$  is \cite{OlverLozierBoisvertEtAl2010} 
[39]
\begin{align}
  W_{\kappa,\mu}(z) & = e^{-\frac{1}{2}z}z^{\frac{1}{2}+\mu} U\Bigl(\tfrac{1}{2}+\mu-\kappa, 1+ 2\mu,z\Bigr)
                  & 2\mu \neq -1, -2, -3, \cdots    
\label{eq:FWequivalent_0070}
\end{align}

By using the following identifications
\begin{align}
   \alpha(s) & =   \tfrac{1}{2}+\mu-\kappa  \\
   \alpha_{+}(s,t) & = \tfrac{1}{2} + \mu
\label{eq:FWequivalent_0071}
\end{align}
such that
\begin{align}
  \kappa & = \alpha_{-}(s,t),
\label{eq:FWequivalent_0072}
\end{align}
and
\begin{align}
   \Psi\bigl(\alpha(s),2\alpha_{+}(s,t), 2\lambda(t-s)\bigr) & = e^{\lambda(t-s)} \bigl(2\lambda(t-s)\bigr)^{-\alpha_{+}(s,t)}
                                                                    W_{\alpha_{-}(s,t),\alpha_{+}(s,t)-\frac{1}{2}}\bigl(2\lambda(t-s)\bigr) . 
\label{eq:FWequivalent_0080}
\end{align}
Hence
\begin{align}                                                                   
 &  \frac{e^{-\lambda(t-s)}(t-s)^{2\alpha_{+}(t,s)-1}}{\Gamma\bigl(\alpha(t)\bigr)}\Psi\bigl(\alpha(s),2\alpha_{+}(s,t), 2\lambda(t-s)\bigr) \nonumber \\
 & \qquad = \frac{e^{-\lambda(t-s)}(t-s)^{2\alpha_{+}(t,s)-1}}{\Gamma\bigl(\alpha(t)\bigr)}
            e^{\lambda(t-s)} \bigl(2\lambda(t-s)\bigr)^{-\alpha_{+}(s,t)} 
            W_{\alpha_{-}(s,t),\alpha_{+}-\frac{1}{2}}\bigl(2\lambda(t-s)\bigr) \nonumber\\
 &  \qquad   = \frac{(t-s)^{\alpha_{+}(t,s)-1}}{\Gamma\bigl(\alpha(t)\bigr)(2\lambda)^{\alpha_{+}(s,t)}}
              W_{\alpha_{-}(s,t),\alpha_{+}-\frac{1}{2}}\bigl(2\lambda(t-s)\bigr)
\label{eq:FWequivalent_0090}
\end{align}
One has \cite{OlverLozierBoisvertEtAl2010}
\begin{align}
  W_{\kappa,\mu}(z) & =  W_{\kappa,-\mu}(z), &  2\mu & \notin \mathbb{Z}. 
\label{eq:FWequivalent_0091}
\end{align}
The condition $2\mu \notin \mathbb{Z}$ requires $\alpha(s) + \alpha(t) \notin \mathbb{Z}$.
By using 
(\ref{eq:FWequivalent_0091}) in (\ref{eq:FWequivalent_0090}),
one thus verify the equality of the covariance functions for the centred Gaussian processes given by 
(\ref{eq:TMfBM_0020}) and (\ref{eq:TMfBM_0040}),
hence the equivalence of the two representations.

\bibliographystyle{unsrtnat}

\bibliography{biblioTFBM}
\end{document}